\documentclass[11pt, twoside]{article}   	
\usepackage{geometry}                		
\usepackage[parfill]{parskip}
\usepackage{graphicx}
\usepackage{amsmath,amsfonts,amssymb}
\usepackage[english]{babel}

\usepackage{algorithm}
\usepackage{algpseudocode}
\usepackage{hyperref}

\newcommand{\R}{\mathbb{R}}
\newcommand{\Lag}{\mathcal{L}}
\DeclareMathOperator{\sech}{sech}

\newcommand{\prt}[2]{\frac{\partial{#1}}{\partial{#2}}}

\DeclareMathOperator*{\argmin}{arg\,min}

\newcommand{\X}{{X}}
\newcommand{\F}{{F}}
\newcommand{\Ftil}{{\tilde{F}}}
\newcommand{\W}{{W}}
\newcommand{\w}{{w}}

\newcommand{\h}{{h}}
\newcommand{\htil}{{\tilde{h}}}
\newcommand{\Win}{{W}_{\rm in}}
\newcommand{\bin}{{b}_{\rm in}}
\newcommand{\HB}{{H}}
\newcommand{\I}{{I}}
\newcommand{\J}{{J}}

\newcommand{\DB}{{D}}
\newcommand{\SB}{{S}}
\newcommand{\U}{{U}}

\newcommand{\G}{{G}}
\newcommand{\Gtil}{{\tilde{G}}}
\newcommand{\BB}{{B}}
\newcommand{\psiB}{{\psi}}
\newcommand{\PsiB}{{\Psi}}
\newcommand{\thB}{{\theta}}
\newcommand{\vthB}{{\vartheta}}
\newcommand{\PhiB}{{\Phi}}
\newcommand{\alphaB}{{\alpha}}
\newcommand{\betaB}{{\beta}}
\newcommand{\GammaB}{{\Gamma}}
\newcommand{\gammaB}{{\gamma}}
\newcommand{\SigmaB}{{\Sigma}}
\newcommand{\e}{{e}}
\newcommand{\Atil}{\tilde{{A}}}
\newcommand{\FI}{Z}
\newcommand{\Bhat}{{\hat{B}}}
\newcommand{\JBhat}{{\hat{J}}}
\newcommand{\GammaBhat}{{\hat{\Gamma}}}
\newcommand{\fhat}{\hat{f}}
\newcommand{\Hhat}{\hat{H}}
\newcommand{\Chat}{\hat{C}}
\newcommand{\Jhat}{\hat{J}}
\newcommand{\FIhat}{\hat{Z}}

\title{Inference of Hamiltonian and Lie-Poisson structures from vector field samples}
\author{Jason E.~Frank\footnote{Utrecht University, Mathematical Institute, Budapestlaan 6, 3584CD Utrecht, the Netherlands, e-mail:\href{mailto:j.e.frank@uu.nl}{j.e.frank@uu.nl}} 
\and Georg A.~Gottwald\footnote{The University of Sydney, School of Mathematics and Statistics, Carslaw Building, Sydney, NSW 2006, Australia, email:\href{mailto:georg.gottwald@sydney.edu.au}{georg.gottwald@sydney.edu.au}\\GAG acknowledges funding from the Australian Research Council, grant DP220100931.}}
\date{\today}							

\begin{document}
\maketitle

\begin{abstract}
We present a method for the data-driven inference of Hamiltonian and Lie-Poisson dynamical systems from vector-field samples, in which both the Hamiltonian and the underlying Poisson tensor are estimated jointly, rather than assuming the structure is known in advance. The Hamiltonian is represented as a linear combination of random feature maps with fixed, randomly drawn internal weights, which renders its estimation a linear regression problem for any given Poisson structure. We treat three cases of increasing generality: a known Poisson tensor, an unknown constant invertible (symplectic) structure matrix, and the case of Lie-Poisson and affine Poisson structures, in which the structure matrix is linear in the state. For the constant case, the Hamiltonian and structure matrix are recovered jointly as the solution of an eigenvalue problem that avoids the trivial zero solution induced by the scaling symmetry of Hamiltonian systems. For Lie-Poisson and affine Poisson systems we use an alternating scheme: an initial Hamiltonian estimate is obtained from a generalized eigenvalue problem for an approximate first integral, and the structure constants are then fit by constrained nonlinear least squares, enforcing the Jacobi identity and, for affine structures, the associated cocycle condition, via continuation in a penalty parameter. Orthogonal matching pursuit selects a compact, well-conditioned subset of random features. Because the learned Hamiltonian decomposes into single-feature terms, the resulting model admits an explicit Poisson/symplectic splitting integrator, which we use to validate the learned systems over long-time integrations via Poincaré sections. We demonstrate the approach on the five-component Lorenz-86 model and a Kirchhoff-type rigid-body system, recovering Hamiltonians, structure constants, and Casimir invariants to high accuracy, and showing convergence of the inferred Poincaré sections to the correct topology of invariant tori as the number of random features increases.
\end{abstract}



\section{Introduction}
                      
Recent years have seen a surge in interest in developing methods to perform data-driven modelling of dynamical systems \cite{LeSchraEtAl23,NoWiSchl23,PiArEtAl2025,PrGe26}. These attempts can broadly be split into two categories, each coming with their own set of applications and aims. One class of methods is concerned with learning surrogate models of a dynamical system that propagate an initial condition to the state at later time. The quest for such surrogate models is driven by computational factors: executing a trained surrogate model offers significant computational advantages over numerically integrating the full high-dimensional often stiff complex dynamical system. A prominent example is weather forecasting where significant advances have recently been achieved by exploiting machine learning \cite{BiEtAl23,LamEtAl23,PriceEtAl24}. 
Recent years have seen advancement in finding surrogate maps for Hamiltonian systems using neural networks, including Lie-Poisson neural networks,  \cite{Eldred24,Putkaradze26} who learn a Poisson map to advance the system, given known Poisson bracket structure, and recent work by Hansen et al.~\cite{HaCeTa25}, who learn a discrete Lagrangian model which they extend to a discrete Lagrangian construction on Lie groups in \cite{HaGhCeMaOw26}.
The second strand of data-driven modelling is concerned with system identification where, given data, one aims to estimate the vector field of the dynamical system.  We are concerned here with the latter class of system identification. Whereas system identification has been successfully achieved for large classes of dissipative dynamics, remarkably even with computationally relatively simple regression methods such as SINDy~\cite{BrPrKu16}, the class of Hamiltonian dynamical systems have proved to be resistant to an efficient estimation. Methods that enjoy success for dissipative systems, including reservoir computers, are unable to preserve the geometric structure innate in Hamiltonian dynamics~\cite{RaMa21,ZhShKa22}. The failure to encode the innate geometrical structure of a Hamiltonian system prevents the trained vector fields being used for long term integration, and typically implies nonphysical energy drain. A promising strategy to overcome these limitations is to approximate the Hamiltonian function rather than the full vector field~\cite{GrDzYo19,HaGlHaMuLa21,DaMe23,BerDiMeKe25,HuOrYi25a,HuOrYi25}. This, however, requires the knowledge of the symplectic structure matrix and excludes the important case of Poisson systems with a state-dependent structure matrix. 
For the particular classes of metriplectic systems or systems with the GENERIC structure neural networks were shown to be successful~\cite{LeTrSt21,ZhShKa22}. Such systems are driven by two integral quantities, a Hamiltonian representing total energy and a single Casimir that is conserved in the dissipation-free limit but grows monotonically in general. 
In \cite{LeTrSt21} fixed structure for the symmetric and antisymmetric tensors is assumed, whereas in \cite{ZhShKa22} the bracket structures may be learned. However, neither article enforces the Jacobi identity, so the skew-symmetric brackets are generally not Poisson. Using Gaussian process kernel regression \cite{HuOrYi25} approximate the Hamiltonians for Poisson systems, given known Poisson structure.

How to estimate simultaneously the Hamiltonian and the Poisson structure remains an open problem. Here we show that this can be achieved for a particular class of Poisson systems, namely, Lie-Poisson and affine Poisson systems for which the structure matrix is linear in the state variables. Hamiltonian system with Lie-Poisson structure \eqref{LPstruct} admit a representation as a dynamical system on the dual to a Lie algebra, and consequently have a rich mathematical theory. They cover a large range of Hamiltonian systems including the free rigid body, ideal incompressible fluids, magnetohydrodynamics, and Vlasov-type kinetic equations \cite{MarsdenRatiu,HolmSchmahStoica}. 

A particularly attractive data-driven modelling framework are random features~\cite{RahimiRecht08}. Random feature maps~\cite{RahimiRecht08} have recently been shown to perform very well for learning surrogate maps for dynamical systems \cite{GottwaldReich21a,NelsenStuart21,LevineStuart22,MaGo25,MaGoNi26}. They are particularly attractive as they render the estimation problem as a regression problem, avoiding costly optimization procedures typical for deep neural networks. Such simple network architectures often outperform more complex machine learning approaches achieving better forecast skill with less hyperparameter tuning and orders of magnitude smaller networks \cite{MaGo25,SchoeWhEtAl26}. Here we use random features to perform system identification. We show that random feature models allow for an effective bilevel regression to estimate both the Hamiltonian and the structure matrix. 

We represent the Hamiltonian as a linear combination of random feature maps. We solve the joint regression for the Hamiltonian and structure constants for symplectic, Lie-Poisson and affine Poisson structures given vector field sample data. 
Ensuring the structure constants adhere to the Jacobi identity places strong constraints on the minimization problem, which we address using a continuation strategy within a nonlinear least squares solver. To test the quality of the learned systems we compute Poincare sections, requiring long numerical integrations. 
Our main conclusions are (1) that representation of the Hamiltonian as a random feature map, in combination with hit-and-run feature selection and optimal matchin pursuit admits a computationally efficient least-squares regression, and (2) combining this approach with constrained regression for the structure constants allows us to infer all of the necessary Hamiltonian structure for accurate system identification and long-term prediction.

The paper is organized as follows. Section~\ref{sec:Ham} recalls the dynamical setting of Hamiltonian systems and in particular of Lie-Poisson systems. In Section~\ref{sec:RFM} we provide a brief overview of random feature maps and develop the random feature map representation  for Hamiltonian systems. In particular, we cover the case when the Poisson structure is known, the case when the Hamiltonian is learned together with a constant symplectic structure, and finally the case of learning both the Hamiltonian and structure constants for the cases of Lie-Poisson and affine Poisson structures. In Section~\ref{sec:splitting} we show that random feature map Hamiltonians admit Poisson splitting methods for their numerical integration. Applying these geometric integrators, we show in Section~\ref{sec:numexp} that our trained random feature models preserve the underlying geometrical structure and are capable of accurately reproducing the dynamics of the four and five-dimensional Lorenz-86 model and a six-dimensional Kirchhoff rigid body model. We conclude in Section~\ref{sec:disc} with a discussion and an outlook.


\section{Hamiltonian dynamical systems and Lie-Poisson systems}
\label{sec:Ham}

System identification is concerned with estimating the vector field $f(x):\R^d \to \R^d$ of a dynamical system 
\begin{equation}\label{ODE}
	\dot{x} = f(x), \qquad x(t) \in \R^d.
\end{equation}
We are concerned here with Hamiltonian systems~\cite{MarsdenRatiu,HolmSchmahStoica}. 
A Hamiltonian system of differential equations on $\R^d$ is defined by a differentiable function $H: \R^d\to\R$ and a symplectic operator $J$, i.e.~a nondegenerate, skew-symmetric operator $J \in \R^{d\times d}: J^\top = -J$:
\begin{equation}\label{HamODE}
	\dot{x} = J \nabla H(x).
\end{equation}
The Poincar\'e theorem establishes that Hamiltonian systems are endowed with the geometric property of symplecticity.  

Symplectic structure places a strong constraint on the dynamics, implying volume preservation, conservation of the Hamiltonian, as well as robustness of integrable structures under small perturbations (KAM theory). Furthermore, excellent symplectic numerical methods are available for simulating Hamiltonian systems, which provide long-time preservation of the above properties~\cite{HaLuWa06,LeRe04,SSCa94}. Any data-driven method estimating the vector field of a Hamiltonian dynamical system hence should preserve the underlying geometric structure and needs to explicitly infer the Hamiltonian structure from sample data. 

Nondegeneracy and skew-symmetry of $J$ together imply that $d$ must be even. In the canonical case that arises naturally in Newtonian mechanics, the structure matrix has the specific form
\begin{equation}\label{Jcanon}
	J = \begin{bmatrix} 0 & I_{d/2} \\ -I_{d/2} & 0 \end{bmatrix},
\end{equation}
where $I_{d/2}$ is the identity matrix of dimension $d/2$.
The phase space coordinate $x$ is partitioned $x = (q,p)$, where $q, p \in \R^{d/2}$ are the conjugate coordinates and momenta, respectively.  On the other hand, if $L \in \R^{d\times d}$ is an invertible matrix, we may define $\tilde{x} = Lx$, $\tilde{H}(\tilde x) = H(L^{-1} \tilde{x})$ and $\tilde{J} = LJL^\top$, and the system \eqref{HamODE} transforms to
\[
	\dot{\tilde{x}} = \tilde{J} \nabla \tilde{H} (\tilde{x}),
\]
where in general $\tilde J$ is a full, skew-symmetric matrix, that typically is not canonical. 

More general Hamiltonian systems may have a Poisson structure, defined by a Poisson bracket, for which the matrix representation yields a system of the form
\begin{equation}\label{PoissonODE}
	\dot{x} = B(x) \nabla H(x),
\end{equation}
where $d$ may be even or odd, and $B(x)$ is a Poisson tensor: a local matrix function of $x$ that must satisfy both the skew-symmetry requirement $B^\top(x) = -B(x)$ and the Jacobi identity:
\begin{equation}\label{jacobiB}
	\sum_\ell \prt{B_{ij}}{x_\ell} B_{\ell k} + \prt{B_{jk}}{x_\ell} B_{\ell i} + \prt{B_{ki}}{x_\ell} B_{\ell j} = 0, \quad \forall x \in \R^d, \forall i,j,k \in \{ 1, \dots, d \}.
\end{equation}
When $B(x)$ is degenerate (if $d$ is odd, for instance), then there may be one or more Casimir functions $Z_i(x)$ satisfying $B(x) \nabla Z_i(x) = 0$. Such functions are first integrals of \eqref{PoissonODE} for any Hamiltonian $H(x)$. 

The Poisson tensor $B(x)$ is characterized by the Darboux theorem, which states that if $B(x)$ has rank $2r\le d$ at $x$, there exists a local change of variables $\tilde{x} = S(x)$ under which $B(x)$ takes the degenerate canonical form
\[
	J = \begin{bmatrix} 0 & I_r & 0 \\ -I_r & 0 & 0 \\ 0 & 0 & 0 \end{bmatrix}.
\] 

A Lie-Poisson structure is defined by a Poisson tensor that is linear in $x$,
\begin{equation}\label{LPstruct}
	B(x)_{ij} = \sum_{k=1}^d C_{ij}^k x_k.
\end{equation}
The skew-symmetry of $B(x)$ translates to the conditions $C_{ij}^k = -C_{ji}^k$ for all $i,j\in \R^d$.
The Jacobi identity \eqref{jacobiB} leads to the conditions
\begin{equation}\label{jacobiIdentity}
 	\sum_{\ell=1}^d C_{ij}^\ell C_{\ell k}^m + C_{jk}^\ell C_{\ell i}^m + C_{ki}^\ell C_{\ell j}^m = 0, \qquad
	i,j,k,m = 1,\dots,d.
\end{equation}
Under a homogeneous linear transformation $\tilde{x} = S(x) = Lx$, with $L$ an invertible matrix, the Lie-Poisson structure constants transform as 
\[
    \tilde{C}_{ij}^k = L_{ii'} L_{jj'} L^{-1}_{k'k} C_{i'j'}^{k'},
\]
which again satisfy \eqref{jacobiIdentity}.

An affine Poisson system is one of the form \eqref{PoissonODE} with a Poisson tensor that is the sum of a Lie-Poisson tensor and a constant skew-symmetric tensor
\begin{equation}\label{affineB}
	B(x)_{ij} = \sum_k C_{ij}^k x_k + J_{ij}.
\end{equation}
In this case, the structure constants $C$ must be skew-symmetric and satisfy \eqref{jacobiIdentity}, and in addition $J$ must be skew-symmetric and satisfy the cocycle condition
\begin{equation}\label{cocycle}
	\sum_{\ell=1}^d C_{ij}^\ell J_{\ell k} + C_{jk}^\ell J_{\ell i} + C_{ki}^\ell J_{\ell j} = 0,
	\qquad \forall i,j,k = 1,\dots,d.
\end{equation}
Affine Poisson structure is again robust under affine linear transformations, i.e.~under $\tilde{x} = S(x) = L x + y$ for invertible $L$.


\section{Random feature maps}
\label{sec:RFM}

Random feature maps can be viewed as a single layer feedforward network in which the internal weights are fixed. Rather than being learned as done in typical feedforward neural networks, in random feature maps the internal weights are drawn from some prescribed random distribution and held fixed. This allows for the estimation of the outer weights by simple least square regression with an explicit analytical expression rather than having to learn them via computationally costly back propagation.

Consider a bounded subset $\mathcal{D}\subset\R^d$ upon which we wish to approximate an unknown function $f: x \mapsto \R^d$ for $x\in\mathcal{D}$; this includes our case of interest where $f(x)$ denotes a vector field $f:\R^d \to \R^d$. We remark that it is easily generalized to functions $f:\R^{n} \to \R^{m}$. Given internal weights $\Win\in\R ^{K\times d}$ and biases $\bin\in\R^K$ of the internal layer drawn from some user-defined distribution, random features $\phi(x)\in \R^K$ are generated by an affine linear transformation of $x\in\R^d$ followed by a nonlinear activation function $\sigma$, that acts component-wise on vectors, as
\begin{align}
\phi(x) = \sigma (\Win x + \bin).
\label{eq:rfm}
\end{align}
We shall use, as often done, $\sigma(z) = \tanh(z)$. The random feature map $\fhat(x)$ approximating a general function $f(x)$ is provided by a linear combination of the $K$ features as
\begin{align}\label{eq:rfmmap}
\fhat:x \mapsto \W\tanh{(\Win x  + \bin)},
\end{align}
where the external layer weights $\W \in\mathbb{R}^{d\times K}$ are to be learned to minimize the data mismatch $||\fhat(x)-f(x)||^2$ over the training data. In particular, we assume we are given $M$ sample data in the sample space $\mathcal{D}$
\[
    X^{(\mu)}\in\mathcal{D}, F^{(\mu)}:\mathcal{D}\to\R^d, \quad \mu=1,\dots,M,
\]
where $F^{(\mu)} = f(X^{(\mu)})$. Introducing the design matrices $\PhiB \in\R^{K \times M}$ and $\F\in \R^{d \times M}$, where the $\mu$th column consists of $\phi(X^{(\mu)})$ and $F^{(\mu)}$, respectively, the estimation of the external weights is given as the minimizer of the following regularized loss function,
\begin{align}
     \underset{\W} {\argmin}\;\|\W\PhiB(\X)-\F\|_F^2 + \beta\|\W\|_F^2
\label{eq:cost},
\end{align}
where $\|\cdot \|_F$ denotes the Frobenius norm. The regularization parameter $\beta$ is a hyperparameter which requires tuning. The solution to the ridge regression problem \eqref{eq:cost} is explicitly given by 
\begin{align}
\W = \F \PhiB^\top \left(\PhiB \PhiB^\top + \beta \I \right)^{-1}.
\label{eq:lsq}
\end{align}
Random feature maps enjoy the universal approximation property, stating that in the limit $K\to\infty$ there exist internal weights $\Win$ and $\bin$ such that they approximate a continuous function arbitrarily well \cite{Cybenko89,ParkSandberg91,Barron93,RahimiRecht08b}. The degree of approximation sensitively depends on the internal weights. If the internal weights are, for example, chosen for the $\tanh$-activation function such that they evaluate to either $\pm 1$ for all $x\in\mathcal{D}$, then the approximation capability is clearly significantly reduced compared to more expressive functions. Recently there has been interest in finding approximate methods to choose the internal parameters to increase the approximation quality of random feature maps given a fixed feature dimension $K$ \cite{LevineStuart22,DunbarEtAl24,MaGoNi26}. We adopt here the strategy proposed in \cite{MaGoNi26} and employ a hit-and-run algorithm to draw the internal parameters such that for given training data the weights do not project the data into either the saturated region of the $\tanh$-function, i.e. evaluating to either $\pm 1$, or into the approximately linear region of the $\tanh$-function. In the former case, the random features cannot differentiate between different data inputs whereas in the latter case they would reduce to a linear model which is incapable of encoding a nonlinear dynamical system. 

In practice the number of random features can be reduced using orthogonal matching pursuit (OMP) as described in \cite{RuZiEl08}. Batch-OMP recursively constructs a subset of features upon which the orthogonal projection of the residual is maximal. OMP terminates when the next best candidate feature is contained in the span of selected features to within a chosen relative tolerance $(10^{-10})$, signalling numerical linear dependence of remaining features.  

In the case when $X$ are samples of state space variables $x$ and $F$ are samples of the associated vector field of a dynamical system $\dot x = f(x)$, one may be tempted to directly apply \eqref{eq:lsq} to obtain an approximation of the vector field $\fhat(x) = \W \tanh(\Win x + \bin)$ with $\W\in\R^{d\times K}$. This strategy, however, would fail for an underlying Hamiltonian system where $f(x) = \BB(x)\nabla H(x)$ in the sense that the particular geometric structure will not be preserved. Rather than expressing the vector field by random features, we follow \cite{HuOrYi25a} and instead express the Hamiltonian $H(x)$ by random features as
\begin{align}
\label{Ham}
	\Hhat(x) = \alphaB^\top \tanh (\Win x + \bin) 
\end{align}
with learnable weights $\alphaB \in \R^K$ and internal weights weights $\Win\in\R ^{K\times d}$ and biases $\bin\in\R^K$. Comparing with the notation of \eqref{eq:rfmmap} we remark that $\W = \alphaB^\top$, i.e.~a row vector for the case of a scalar valued function. 

By expressing the Hamiltonian as a random feature map~\eqref{Ham}, the vector field can be expressed as a linear function of the external weights $\alpha_\kappa$. 
Define 
\begin{align}
\psi_\kappa(x) = \sech^2 (\w_\kappa^\top x + b_\kappa), \qquad \kappa=1,\dots,K,
\label{eq:psik}
\end{align}
where $\w_\kappa^\top \in \R^d$ denotes the $\kappa$th row of $\Win\in\R ^{K\times d}$ and $b_\kappa$ is the $\kappa$th component of the biases $\bin\in\R^K$. 

Introduce the associated design matrix $\PsiB \in\R^{K\times M}$ with components $\psi_\kappa^{(\mu)} = \psi_\kappa (X^{(\mu)})$, $k=1,\dots,K$, $\mu=1,\dots,M$.
Let $\psiB^{(\mu)}=\psiB(X^{(\mu)})\in\R^K$ denote the $\mu$th column of $\PsiB$. The gradient of the Hamiltonian expressed as a random feature map with $\phi_\kappa(x)$~\eqref{eq:rfm} with $\sigma=\tanh$, is then given by
\begin{align}
\nabla \Hhat (X^{(\mu)}) 
= \sum_{\kappa=1}^K \alpha_\kappa \psi_\kappa^{(\mu)} \w_\kappa 
= \Win^\top {\rm{diag}}(\alphaB) \psiB^{(\mu)}.
\end{align}

In the following we describe how to use random feature maps to learn three different classes of Hamiltonian dynamical systems of increasing complexity: (i) a generic Hamiltonian system with a given Poisson tensor $\BB(x)$, for which we only need to estimate the Hamiltonian, (ii) an unknown symplectic (i.e.~constant, nondegenerate, and skew-symmetric) structure matrix for which we need to estimate the Hamiltonian and the structure matrix, and (iii) an unknown general affine or Lie-Poisson system for which we estimate the Hamiltonian and structure constants.


\subsection{Learning a Hamiltonian system with known Poisson structure}
\label{sec:Honly}

For a Hamiltonian system with a known Poisson tensor $B(x)$ the loss function~\eqref{eq:cost} is given as a quadratic equation for the unknown coefficients $\alphaB \in \R^K$ determining the unknown Hamiltonian~\eqref{Ham}. Introducing the Poisson tensor evaluated at sample points $\BB^{(\mu)} = \BB( X^{(\mu)})$ we solve the following minimization problem 
\begin{align}
     \underset{\alphaB}{\argmin}\;
     \sum_{\mu=1}^M \|\BB^{(\mu)} \Win^\top {\rm{diag}}(\alphaB) \psi^{(\mu)} - \F^{(\mu)}\|_F^2 + \beta\|\alphaB\|_2^2
\label{eq:cost_canonHam}.
\end{align}
The regularization term penalizing large norms of $\alphaB$ serves two related purposes: first it renders a potentially ill-conditioned problem solvable, and second it prevents the occurrence of large values of entries $\alpha_\kappa$ which would cause the approximate Hamiltonian vector field to be stiff, requiring undesirable small time steps for its integration. 

Introducing the ridge-regularized Gram matrix
\begin{align}
\G = \sum_{\mu =1}^M  ( \Win {\BB^{(\mu)}}^\top \BB^{(\mu)} \Win^\top) \odot (\psiB^{(\mu)} {\psiB^{(\mu)}}^\top) + \beta \I_{K}\in \R^{K\times K},
\label{eq:G}
\end{align}
where $\odot$ denotes the Hadamard product, and
\begin{align}
\h = \sum_{\mu=1}^M (\Win{\BB^{(\mu)}}^\top \F^{(\mu)} )\odot \psiB^{(\mu)} \in \R^K,
\label{eq:h}
\end{align}
the normal equation with respect to $\alphaB$ of the minimization problem~\eqref{eq:cost_canonHam} is given by
\begin{align}
    \G\alphaB = \h,
    \label{eq:normEq_alpha}
\end{align}
with solution given by 
\begin{align}\label{alphaSol}
    \alphaB = \G^\dagger \h,
\end{align}
where $\G^\dagger$ denotes the pseudoinverse of $\G$. 

In the case when $\BB(x)=\J$ is constant, we obtain
\begin{align}
\G = ( \Win\J^\top\J\Win^\top) \odot (\PsiB \PsiB^\top) + \beta \I_{K}\in \R^{K\times K},
\label{eq:G_Bconst}
\end{align}
where we used that $\sum_{\mu=1}^M\psiB^{(\mu)}{\psiB^{(\mu)}}^\top=\Psi\Psi^\top$, and
\begin{align}
\h = {\rm{diag}}(\Win\J^\top \F \PsiB^\top) \in \R^K.
\label{eq:h_Bconst}
\end{align}
This includes the case of a canonical Hamiltonian system, $\J$ given by \eqref{Jcanon}.

In our numerical implementation we do not compute the normal equation \eqref{eq:normEq_alpha}, but instead work directly with the overconstrained equation for $\alphaB$, computing the singular value decomposition of the associated design matrix $\PsiB$ and including the ridge parameter in the pseudoinverse. The solution is mathematically equivalent to \eqref{alphaSol}.


\subsection{Learning a Hamiltonian system with unknown constant nondegenerate structure matrix}
\label{sec:nondegenHam}

If our goal is to infer a general Hamiltonian system from vector field samples, it is unlikely that our measurements will be conveniently grouped into canonically conjugate pairs. Under the assumption that the symplectic structure is nondegenerate (which requires $d$ to be even), there exists a linear transformation to canonical form and the approach of the previous section may be used. Alternatively, we may infer a symplectic structure matrix $\J$ with $\GammaB = \J^{-1}$ together with the Hamiltonian by solving the least squares problem 
\begin{align}
     \underset{\alphaB,\GammaBhat}{\argmin}\;\|\Win^\top {\rm{diag}}(\alphaB) \PsiB - \GammaBhat \F\|_F^2 + \beta\|\alphaB\|_2^2
\label{eq:cost_symplectic}.
\end{align}
under the constraint that $\GammaBhat$ is skew-symmetric. The gradient of the loss function is linear in both $\alphaB$ and $\GammaBhat$. The trivial minimizer $\GammaBhat = 0$, $\alphaB = 0$ reflects the scaling symmetry of the Hamiltonian system: $H \mapsto \xi H$, $J \mapsto \xi^{-1} J$. To obtain a unique solution we require without loss of generality that $\|\GammaBhat \|_F^2=1$.

Since $\GammaB = -\GammaB^\top$ there are only $q=\binom{d}{2}$ distinct entries of $ \GammaB$. Let $r = r(i,j)$, $i<j$ denote an ordering, choose a basis $\SB_1,\ldots,\SB_q$ of the space of skew-symmetric matrices with
\begin{align}
\SB_{r(i,j)} = \frac{1}{\sqrt{2}}(\e_i\e_j^\top - \e_j\e_i^\top) \qquad i<j,
\end{align}
and write $\GammaBhat = \sum_{r=1}^q \gammaB_r \SB_r$ with $\gammaB\in\R^q$. Note the $\SB_r$ are scaled such that $\|\gammaB\|_2^2=1$ implies $\|\GammaB\|_F^2=1$. Substituting into the data mismatch term of the loss function~\eqref{eq:cost_symplectic}, we obtain
\begin{align}
\|\Win^\top {\rm{diag}}(\alphaB) \PsiB - \GammaBhat \F\|_F^2 = \|\sum_{\kappa=1}^K \alpha_\kappa \w_\kappa\psiB_k^\top - \sum_{r=1}^q \gammaB_r\SB_r\F \|_F^2.
\end{align}

The normal equation for $\alphaB$, analogous to \eqref{eq:normEq_alpha}, is
\begin{equation}
    \Gtil_\beta \alphaB = \htil,
\end{equation}
where the Gram matrix in this case is given by 
\begin{align}
\Gtil_\beta = \sum_{\mu =1}^M  ( \Win \Win^\top) \odot (\psiB^{(\mu)} {\psiB^{(\mu)}}^\top) + \beta \I_{K}\in \R^{K\times K}.
\label{eq:Ghat}
\end{align}

Since $\htil$ is linear in $\GammaBhat$ (cf~\eqref{eq:h}), it is linear in $\gammaB$. We find
\begin{align}
\htil(\GammaB) = \DB \gammaB ,    
\end{align}
for $\DB\in \R^{K\times q}$ with 
\begin{align}
D_{k r} = \w_k^\top \SB_r\F \psiB_k,
\end{align}
where $\psiB_k \in \R^M$ is the transpose of the $k$-th row of $\PsiB$. 

Since $\DB \gammaB$ is in the range of $\Gtil_\beta$ we can write the minimizer of \eqref{eq:cost_symplectic} with respect to $\alphaB$ as 
\begin{equation}
\alphaB = \Gtil_\beta^\dagger \DB \gammaB.
\end{equation}

Since $\alphaB^\top \Gtil_\beta \alphaB = \alphaB^\top \DB \gammaB$ for the minimizer $\alphaB$, the data mismatch loss function, evaluated at the minimizer, becomes
\begin{align}
\mathcal{L}(\gammaB)  
&= \|\Win^\top {\rm{diag}}(\alphaB) \PsiB - \GammaBhat \F\|_F^2 \nonumber \\
&= \gammaB^\top \HB \gammaB - 2 \alphaB^\top \DB \gammaB + \alphaB^\top \Gtil_\beta \alphaB \nonumber \\
&= \gammaB^\top \HB \gammaB - \alphaB^\top \DB \gammaB \nonumber \\
& = \gammaB^\top (\HB - \DB^\top\Gtil_\beta^\dagger \DB) \gammaB,
\end{align}
where we introduced $\HB\in \R^{q\times q}$ with
\begin{align}
\HB_{rs} = {\rm{tr}}(\SB_r\F\F^\top\SB_s^\top), \end{align}
and 
where we used 
that ${\Gtil}_\beta^\dagger$ is symmetric. 

Hence the minimum is given as the eigenvector corresponding to the smallest eigenvalue $\lambda$ of
\begin{align}\label{Jeig}
    (\HB-\DB^\top \Gtil_\beta^\dagger \DB)\gammaB = \lambda \gammaB,
\end{align}
with $\|\gammaB\|_2^2=1$.

In our numerical implementation we solve \eqref{Jeig} for $\gamma$, construct $\JBhat = \GammaBhat^{-1}$ and rescale such that $\max_{ij} |\Jhat_{ij}| = 1$.  Thereafter we apply OMP to reduce the feature pool. Subsequently we solve for $\alphaB$ using the methodology of Section \ref{sec:Honly}. 


\subsection{Learning a Hamiltonian system with Lie-Poisson or affine Poisson structure \label{sec:LPstruct}}

In coordinates in $\R^d$ the Lie-Poisson structure can be considered a special case ($\J=0$) of the affine Poisson structure \eqref{affineB}. In this section we consider the inference of the more general affine Poisson structure \eqref{affineB} along with the Hamiltonian. 

The unknown Poisson tensor $\BB(x)$ defined in \eqref{affineB} is linear in the structure coefficients $C_{ij}^k$ and $J_{ij}$, which in turn are subject to the constraints imposed by the Jacobi identity \eqref{jacobiIdentity} and the co-cycle condition \eqref{cocycle}. We parameterize $\Bhat$ by its distinct coefficients up to skew-symmetry, denoting $\thB = (\Chat_{ij}^k, \Jhat_{ij})$.
We define the loss function 
\begin{equation}\label{L3}
 	\Lag(\alphaB,\thB) = \sum_{\mu=1}^M \| \Bhat(X^{(\mu)}) \nabla \Hhat (X^{(\mu)}) - F^{(\mu)} \|_F^2, 
\end{equation}
where $\Hhat$ is linear in the coefficients $\alphaB$, and $\Bhat$ is linear in $\thB$ subject to constraints $g(\thB)=0$ that encode \eqref{jacobiIdentity} and \eqref{cocycle}.

To solve for $\alphaB$ and $\thB$ we iterate, constructing the sequence $(\thB^{(n)}, \alphaB^{(n)})$ by first minimizing $\mathcal{L}$ with respect to $\thB$ for fixed $\alphaB$ subject to \eqref{jacobiIdentity} and \eqref{cocycle}, then solving for $\alphaB$ by minimizing \eqref{eq:cost_canonHam} with explicit solution~\eqref{alphaSol} as described in Section \ref{sec:Honly}. This process of bilevel optimization is repeated until convergence. To start the iteration we require an initial guess $\alphaB^{(0)}$, which we construct using the approach of the following subsection.


\subsubsection{Initialisation step: Estimating a first integral from vector field data}

We adapt the methodology from \cite{ZhZhKe23}. A smooth first integral $\FI(x):\R^d\to\R$ of \eqref{ODE} satisfies $\nabla \FI(x)^\top f(x) = 0$ for all $x$.  However, any such function is nonunique, and if there are multiple first integrals $\FI_1$, \dots $\FI_r$, then $q(\FI_1,\dots,\FI_r)$ is also a first integral for any smooth $q$. Moreover, $\FI(x)={\it{const}}$ is trivially a constant of motion. 

We normalize the vector field data, defining $\Ftil^{(\mu)} = F^{(\mu)}/\| F^{(\mu)} \|$.
To find a non-trivial first integral we express it as a linear combination of random features analogously to~\eqref{Ham} with coefficients $\alphaB^{(0)}$, and we seek $\alphaB^{(0)}$ minimizing the constrained loss function
\begin{align}
	\Lag_{\rm{init}}(\alphaB,\lambda) = \sum_{\mu=1}^M \| \nabla \FIhat(X^{(\mu)})^\top  \Ftil^{(\mu)} \|_2^2 - \lambda \left( \sum_{\mu=1}^M (\| \nabla \FIhat( X^{(\mu)}) \|_2^2)-1 \right),
\end{align}
where $\FIhat$ is linear in $\alphaB^{(0)}$ and $\lambda$ is a Lagrange multiplier for the constraint, which precludes trivial constant first integrals with $\nabla \FIhat(x)\equiv 0$. We rewrite the loss function explicitly as 
\begin{align}
    \Lag_{\rm{init}}(\alphaB,\lambda) = \sum_{\mu=1}^M  ||\Ftil^{(\mu)\top}\Win^\top {\rm{diag}}(\psiB^{(\mu)})\alphaB||_F^2 - \lambda (\sum_{\mu=1}^M ||\Win^\top  {\rm{diag}}(\psiB^{(\mu)})\alphaB||_2^2-1).
\end{align}
Introducing
\begin{align}
    \Atil &= \sum_{\mu=1}^M  {\rm{diag}}(\psiB^{(\mu)})\Win \Ftil^{(\mu)} \Ftil^{(\mu)\top} \Win^\top  {\rm{diag}}(\psiB^{(\mu)})\nonumber \\
    &= [ \PsiB \odot (\Win\Ftil)][ \PsiB \odot(\Win \Ftil)]^\top \in \R^{K\times K},
\end{align}
the loss function reads
\begin{align}
    \Lag_{\rm{init}}(\alphaB;\lambda) = \alphaB^\top \Atil \alphaB - \lambda(\alphaB^\top \Gtil_0 \alphaB -1),
\end{align}
where $\Gtil_0$ is defined by \eqref{eq:Ghat} with $\beta=0$.
The normal equation with respect to $\alphaB^{(0)}$ is given by
\begin{align}
    \Atil \alphaB^{(0)} = \lambda \Gtil_0 \alphaB^{(0)},
    \label{eq:alpha_genEV}
\end{align}
and minimization of $\Lag_{\rm{init}}$ amounts to finding the generalized eigenvector corresponding to the smallest generalized eigenvalue $\lambda$.

In general $\Gtil_0$ is not invertible. We decompose $\Gtil_0 = \U_r \SigmaB_r\U_r^\top$ such that $\SigmaB_r$ contains all $r\le K$ nonzero eigenvalues and $\U_r\in \R^{K\times r}$ contains the first $r$ eigenvectors of $\Gtil_0$, i.e. $\U_r^\top \U_r = \I_r$. We solve for $\alphaB^{(0)}$ in the column space of $\U_r$. 
Substituting $\alphaB^{(0)} = \U_r \betaB$ for $\betaB\in \R^r$, the generalized eigenvalue problem~\eqref{eq:alpha_genEV} becomes
\begin{align}
    \Atil \U_r\betaB = \lambda \U_r\SigmaB_r\betaB.
\end{align}
Introducing $\tilde \betaB = {\SigmaB_r}^{\tfrac12} \betaB$ and symmetrizing, we obtain 
\begin{align}
    \SigmaB_r^{-\tfrac12} \U_r^\top \Atil \U_r \SigmaB_r^{-\tfrac12} \tilde \betaB = \lambda \tilde \betaB.
\end{align}
This allows us to stably solve for the eigenvector $\tilde \betaB_{\rm{min}}$ corresponding to the smallest eigenvalue $\lambda=\lambda_{\rm{min}}$ of a symmetric matrix. Transforming back into original coordinates, we finally obtain
\begin{equation}\label{eqn:alphaguess}
    \alphaB^{(0)} = \U_r \SigmaB_r^{-\tfrac12} \tilde \betaB_{\rm{min}},
\end{equation}
which we use to construct the initial guess for the Hamiltonian~\eqref{Ham}. 

This procedure to find an initial guess of the Hamiltonian finds an approximation to a first integral that in general will be a function of the Hamiltonian as well as any Casimirs and or Noether momenta present. Whereas Casimirs are irrelevant since they vanish in the Poisson bracket, the presence of Noether momenta could conceivably degrade 
the quality of the first guess. Nevertheless, the iteration of the following section may still converge.


\subsubsection{Nonlinear least squares for structure constants\label{sec:nllsqr}}

In the $n$th outer iteration of our bilevel iteration process, given an estimate for $\alphaB^{(n-1)}$ we minimize \eqref{L3} for the structure constants 
$\thB=(\Chat_{ij}^k,\Jhat_{ij})$ of $\Bhat(x)$, subject to constraints imposed by the Jacobi identity
\eqref{jacobiIdentity} and co-cycle condition \eqref{cocycle}, encoded in the condition $g(\thB) = 0$.
In this section we propose a method for computing $\thB^{(n)}$ given $\alphaB^{(n-1)}$.

We emphasize that this problem is challenging. The number of distinct (up to skew-symmetry) coefficients $C_{ij}^k$ is $d \cdot \binom{d}{2}$ and of $J_{ij}$ is $\binom{d}{2}$, i.e.,~the combined set of parameters $\thB$ grows as $\mathcal{O}(d^3)$. 
The number of distinct constraints imposed by the Jacobi identity \eqref{jacobiIdentity} is $d \cdot \binom{d}{3}$ and by the cocycle condition \eqref{cocycle} is $\binom{d}{3}$, i.e.,~the number of constraints grows as $\mathcal{O}(d^4)$. The number of constraints precisely equals the number of parameters for $d=5$.
Consequently, for $d>5$ an increasingly large number of constraints must converge when we find a solution. 
The algebraic variety $g(\thB)=0$ defined by the constraints is reducible and stratified with singular connections. This nonsmoothness presents an obstacle to methods that rely on gradient minimization with projection onto the variety. Instead we employ a weak constraint with a continutation strategy in the constraint parameter $\lambda$. Starting from a point off the variety that simply minimizes the data mismatch, we approach the variety by gradually increasing $\lambda$ until the Jacobi identity and cocycle conditions are satisfied.

For fixed $\alphaB=\alphaB^{(n-1)}$ we formulate \eqref{L3} abstractly as a weak-constraint nonlinear least squares problem for $\thB$
\begin{equation}\label{L3pp}
	\argmin_{\thB} \Lag(\thB,\alphaB^{(n-1)}) + \lambda \| g(\thB) \|_2^2.
\end{equation}
where $\lambda$ is a weak constraint parameter.  

We solve \eqref{L3pp} using continuation over  $\lambda$, where for small $\lambda$ we fit the data well and as $\lambda$ increases we enforce the constraint.
Specifically, we solve \eqref{L3pp} using an inner iteration over a schedule of $\lambda$ values of increasing magnitude $\lambda_j = 10^{3j-4}$, $j=1,\dots,j_{\text{max}}$. 

For fixed $\lambda_j$, $j>0$ we compute $\vthB^{(j)}$ from \eqref{L3pp} using Levenberg-Marquardt iterations applied to the residual 
\[
    r(\thB) = D\Lag(\thB) + \lambda g(\thB) Dg(\thB) = 0
\]
from a starting value $\vthB^{(j-1)}$, the output from the previous iteration computed with $\lambda_{j-1}$. 
The Levenberg-Marquardt iteration requires the Hessian of $g$, which is a constant matrix.  

Let $\vthB = h_\mathrm{LM}(\vthB_0, \lambda; \alphaB)$ denote the solution of the Levenberg-Marquardt iteration with starting value $\vthB_0$, weak constraint parameter $\lambda$, and fixed $\alphaB$. Then we solve the schedule
\[
	\vthB^{(j)} = h_\mathrm{LM}(\vthB^{(j-1)},\lambda_j; \alphaB^{(n-1)}), \qquad j=1,\dots,j_{\text{max}},
\]
with $\vthB^{(0)}$ given by solving the unconstrained (linear) least squares problem for the parameters. Having completed the full schedule, we define $\thB^{(n)} = \vthB^{(j_{\text{max}})}$.
In all numerical experiments in the following Sections we set $j_{\text{max}}=8$. 

Given $\thB^{(n)}$ we compute $\alphaB^{(n)}$ by minimizing \eqref{eq:normEq_alpha} using the methodology of Section \ref{sec:Honly}.

We summarize the algorithm in the following pseudocode:

\begin{algorithm}[!h]
\caption{Method for jointly learning Hamiltonian and affine Poisson structure constants $\thB=(J_{ij},C_{ij}^k)$}
\begin{algorithmic}[1]
\State Solve \eqref{eqn:alphaguess} for $\alphaB^{(0)}$
\For{$n=1,\dots$}
    \State $\vthB^{(0)} \gets \arg\min_{\thB} \Lag_3(\thB)$ from  \eqref{L3pp} with $\lambda=0$.
    \For{$j=1,\dots,j_{\mathrm{max}}$}
        \State $\lambda_j \gets 10^{3j-4}$
        \State $\vthB^{(j)} \gets h_\mathrm{LM}(\vthB^{(j-1)},\lambda_j; \alphaB^{(n-1)})$
    \EndFor
    \State $\thB^{(n)} \gets \vthB^{(j_{\mathrm{max}})}$
    \State Solve \eqref{eq:normEq_alpha} for $\alphaB^{(n)}$
\EndFor
\end{algorithmic}
\end{algorithm}


\section{Numerical integration of the inferred model via a Poisson splitting}
\label{sec:splitting}

Having inferred a Hamiltonian or Lie-Poisson vector field from sample data, one may wish to make predictions by numerically integrating the inferred system. It is generally advised to use numerical methods that are themselves symplectic or Poisson maps when integrating such systems \cite{SSCa94,HaLuWa06,LeRe04}. One class of numerical schemes that are symplectic for general constant, skew-symmetric Poisson tensors is the Gauss-Legendre family of Runge-Kutta methods. However, they have the drawback of being implicit, requiring computationally costly root solvers. 
Alternatively, splitting methods lead to explicitly computable methods that are applicable to both symplectic and Poisson flows. 

We begin with symplectic Hamiltonian systems. Because the exact flow of a Hamiltonian system is a symplectic map (with respect to the particular structure matrix $\J$ (or its learned approximation $\JBhat$)), and the composition of symplectic maps is again a symplectic map, we split the Hamiltonian as a sum of individual random feature maps:
\begin{equation}\label{splitHam}
	\Hhat(x) = \sum_{\kappa=1}^K  \Hhat_\kappa(x), \qquad \Hhat_\kappa(x) = \alpha_\kappa \tanh (\w_\kappa^\top x + b_\kappa).
\end{equation}
The flow generated by the lone feature $\Hhat_\kappa$ is the solution of 
\begin{equation}\label{splitJ}
	\dot{x} = \J \nabla \Hhat_\kappa(x) = \alpha_\kappa \sech^2(\w_\kappa^\top x + b_\kappa) \J \w_\kappa
\end{equation}
Due to the skew symmetry we see that
\[
	\w_\kappa^\top \dot{x} = \alpha_\kappa\sech^2(\w_\kappa^\top x + b_\kappa) \w_\kappa^\top \J \w_\kappa = 0.
\]
Consequently, $\w_\kappa^\top x$ is a constant of motion of this flow, and the right side of \eqref{splitJ} is constant. Therefore, the (exact) flow map associated with each Hamiltonian $\Hhat_\kappa$ is simply given by
\[
	\varphi_t^{(\kappa)} (x_0) = x_0 + t \alpha_\kappa\sech^2(\w_\kappa^\top x_0 + b_\kappa) \J \w_\kappa .
\]
The Lie-Trotter composition 
\begin{equation} \label{Trotter}
	\phi_{\Delta t} (x) = \varphi_{\Delta t}^{(1)} \circ \cdots \circ \varphi_{\Delta t}^{(K)} (x)
 \end{equation}
defines a first order symplectic numerical integrator. Its adjoint is the composition in  reversed order
\[
 	\phi_{\Delta t}^* (x) = \varphi_{\Delta t}^{(K)} \circ \cdots \circ \varphi_{\Delta t}^{(1)} (x).
\]
For our numerical computations we use the symmetric Strang composition
\begin{equation}\label{Strang}
 	\Phi_{\Delta t}(x) = \phi_{\Delta t/2} \circ \phi_{\Delta t/2}^* (x),
\end{equation}
which is symplectic and second order accurate \cite{HaLuWa06}.

Note that the cost of the first order splitting \eqref{Trotter} is essentially the same as the evaluation of the whole vector field, and the symmetric splitting is twice as expensive. Consequently, these methods have comparable costs to first and second order explicit Runge-Kutta methods applied to the learned system.

This splitting scheme generalizes recursively to larger $K$.  For large $K$ the splitting errors may become significant, typically requiring a smaller stepsize than is needed to integrate the true Hamiltonian in some examples (see below).

For inferred affine and Lie-Poisson systems, the arguments are similar. The $\kappa$th term in \eqref{splitHam} generates the dynamics
\[
	\dot{x} = \BB(x) \nabla \Hhat_\kappa(x) = \alpha_\kappa\sech^2(\w_\kappa^\top x + b_\kappa) \BB(x) \w_\kappa.
\]
Again we find $\w_\kappa^\top \dot{x} = 0$ by skew-symmetry. Recalling the definition~\eqref{eq:psik} with $\psi_\kappa(x) = \sech^2(\w_\kappa^\top x + b_\kappa)$ and denoting by $Q_\kappa = C(\cdot,\w_\kappa,\cdot)$ the matrix obtained by contracting the structure tensor with $w_\kappa$ along its second index and $y_\kappa = \J \w_\kappa$, we can write
\[
	\dot{x}  = \alpha_\kappa\psi_\kappa(x)(Q_\kappa x + y_\kappa),
\]
where $\psi_\kappa(x)$ is constant along the flow.  Consequently, the flow generated by $\Hhat_\kappa$ is given by the variation of constants formula, formally as
\begin{equation}\label{LPsubstep}
	\varphi_t^{(\kappa)}(x_0) = e^{t \alpha_\kappa\psi_\kappa(x_0) Q_\kappa}x_0 + Q_\kappa^{-1} (I - e^{t \alpha_\kappa\psi_\kappa(x_0) Q_\kappa}) y_\kappa.
\end{equation}
which can be implemented by computing a matrix exponential (see \cite{HaLuWa06} for efficient implementation, invertibility of $Q_\kappa$ is not required).
The resulting Lie-Trotter \eqref{Trotter} and Strang \eqref{Strang} splittings, defined identically, yield Poisson integrators of first and second order, respectively. Of course, computing $2K$ matrix exponentials per time step may become costly.


\section{Numerical results}
\label{sec:numexp}
We shall now show that the methodology of system identification using random feature maps is capable of reliably inferring Hamiltonian dynamical systems and their underlying symplectic, Lie-Poisson or affine Poisson structure. We consider two models: The Lorenz 86 model which comes in two flavours; a $5$-dimensional version with a nontrivial affine Poisson structure and a $4$-dimensional symplectic version. We further consider the $6$-dimensional Kirchhoff equations which are a Lie-Poisson system. 

We will employ the methods described in the previous sections to learn the Hamiltonian and/or the underlying structure matrix. In all computations we assume noise-free vector field data. The internal weights for the random feature maps, $\Win$ and $\bin$ are sampled using the hit-and-run algorithm~\cite{MaGoNi26}.


\subsection{The models}
We first introduce the models and define the Poincar\'e section which will be used to assess the accuracy of the learned models. 


\subsubsection{The Lorenz 1986 model}
\label{sec:mod_l86}

In \cite{Lorenz86} E.N.~Lorenz constructed a low-dimensional model based on the rotating shallow water equations to study simplified dynamics of geostrophic balance. The model equations (rescaled after \cite{BoSh96}) are
\begin{subequations}\label{L86}
\begin{align}
	\dot{x}_1 &= -x_2 x_3 + b x_2 x_5,\\
	\dot{x}_2 &= x_1x_3 - b x_1x_5,\\
	\dot{x}_3 &= -x_1 x_2, \\
	\dot{x}_4 &= - \frac{1}{\varepsilon} x_5,\\
	\dot{x}_5 &= \frac{1}{\varepsilon} x_4 + b x_1 x_2,
\end{align}
\end{subequations}
where $b$ represents the rotational Froude number and $\varepsilon = b\rho(1 + b^2)^{1/2}$ is a perturbation parameter (here $\rho$ is the Rossby number; we assume $\rho=1$). The system \eqref{L86}, in the following referred to as the L86 system, is a Poisson system with affine Poisson structure \eqref{affineB}, with nonzero elements
\begin{equation} \label{L86affineB}
	\{ C_{13}^2 = -1, C_{15}^2 = b, C_{23}^1 = 1, C_{25}^1 = -b, \}, \qquad  J_{45} = -1/\varepsilon,
\end{equation}
up to skew-symmetry, and with Hamiltonian
\begin{equation}\label{L86ham}
	H(x) = \frac{1}{2} ( x_1^2 + 2x_2^2 + x_3^2 + x_4^2 + x_5^2).
\end{equation}
In addition the L86 system~\eqref{L86} preserves the Casimir
\begin{equation}\label{L86cas}
	Z(x) = \frac{1}{2} (x_1^2 + x_2^2).
\end{equation}

By restricting to a level set of constant Casimir $Z$, Bokhove and Shepherd \cite{BoSh96} 
expressed the $5$-dimensional L86 model as a canonical symplectic $4$-dimensional model given by 
\begin{subequations}\label{L86canon}
\begin{align}
	\dot{q}_1 &= \varepsilon p_1,\\
	\dot{q}_2 &= p_2, \\
	\dot{p}_1 &= -\varepsilon Z \sin 2(q_1 + b\sqrt{\varepsilon}q_2), \\
	\dot{p}_2 &= -q_2 - \varepsilon^{3/2}bZ  \sin 2(q_1 + b\sqrt{\varepsilon}q_2).
\end{align}
\end{subequations}
The canonical coordinates are defined with respect to Lorenz's formulation by
\[
	q_1 = \phi - \varepsilon b x_4, \quad p_1 = x_3, \quad q_2 = \sqrt{\varepsilon} x_4, \quad p_2 = -\sqrt{\varepsilon} x_5,
\]
where $\phi$ is the phase with respect to the Casimir $Z$, i.e.~$x_1 = \sqrt{2Z} \cos \phi$, $x_2 = \sqrt{2Z} \sin \phi$.  The canonical Hamiltonian of the reduced system is given by 
\begin{equation}\label{Hcanon}
	H_{\mathrm{canon}}(q_1,q_2,p_1,p_2) = \varepsilon H = \frac{\varepsilon}{2} \left[ p_1^2 - Z \cos 2(q_1 + b\sqrt{\varepsilon} q_2) \right] + \frac{1}{2} (p_2^2 + q_2^2).
\end{equation}
We will use this $4$-dimensional representation of the L86 model in the following. 

As a test problem for the methods, we compute a Poincar\'{e} section for the parameter choice 
\[
	H = b = \frac{1}{2}, \qquad \rho = Z = 1,
\]
which supports both invariant tori and chaotic dynamics. We choose a Poincar\'{e} section corresponding to the $(\phi,x_3)$ plane at $x_4 = 0$ and $\dot{x}_4>0$  (in canonical coordinates, the $(q_1,p_1)$ plane at $q_2=0$, $p_2>0$).  

Computing the Poincar\'{e} section is a challenging test problem because it requires a sufficiently accurate approximation of the Hamiltonian and symplectic/Poisson structure to allow for an accurate long time integration. Significant errors are immediately visible in the form of distortions or topological changes to the invariant tori. 

\begin{figure}[!ht]
\centering
\includegraphics[width=0.5\textwidth]{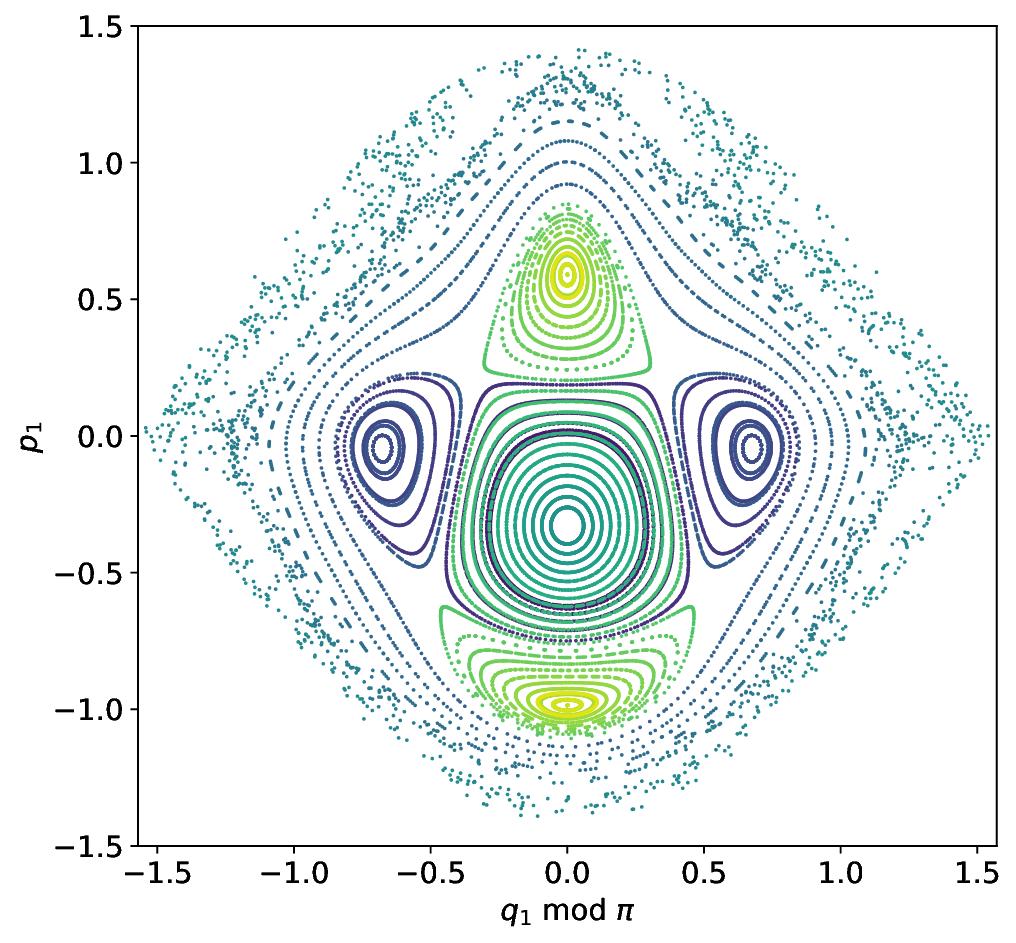}
\caption{Poincar\'{e} section of the exact L86 model \eqref{L86canon} computed with the St\"{o}rmer-Verlet method using a stepsize $\Delta t = 0.0025$.}
\label{fig:exactPoincare}
\end{figure}

The Poincar\'{e} sections are computed using simulations to time $T=2000$ (canonical system time scale) from a set of 50 initial conditions along the two rays:
\[
	p_1 = 0, \quad q_1 \in \{ j \Delta q, j=1,\dots,25 \}, \quad \Delta q = \frac{\pi}{54},
\]
and 
\[
	q_1 = 0, \qquad p_1 \in \{ -0.3 + j \Delta p, j=1,\dots,25 \}, \quad \Delta p = \frac{0.3 + 1/\sqrt{2}}{27}.
\]	
Figure \ref{fig:exactPoincare} shows the Poincare section computed in canonical coordinates using the St\"ormer-Verlet method for the canonical L86 system~\eqref{L86canon}.

\subsubsection{The Kirchhoff equations}
\label{sec:kirchhoff}

We further consider a Kirchhoff-type  model for rigid body motion in an ideal fluid, defined on $\R^6$. 
\begin{equation} \label{kirchhoff}
    \dot{\Omega} = -\Omega \times \prt{H}{\Omega} -\Gamma \times \prt{H}{\Gamma}, \qquad \dot{\Gamma} = - \Gamma \times \prt{H}{\Omega},
\end{equation}
where $\Omega(t), \Gamma(t) \in \R^3$ are the angular and linear velocities, and the Hamiltonian is defined by
\[
    H(\Omega,\Gamma) = \frac{1}{2} \sum_{i=1}^3 a_i \Omega_i^2 + \sum_{i=1}^3 b_i \Omega_i \Gamma_i + \frac{1}{2} \sum_{i=1}^3 d_i \Gamma_i^2
\]
The system is Lie-Poisson and admits two Casimir functions
\begin{equation}\label{KirchCasimirs}
    Z_1 = \| \Gamma \|^2, \qquad Z_2 = \Omega\cdot\Gamma.
\end{equation}


\subsection{Inferring the Hamiltonian for a given canonical structure}

We first study inference of the Hamiltonian assuming that the canonical symplectic structure is known.  We consider here the canonical L86 model~\eqref{L86canon}. Given vector field data we learn the canonical Hamiltonian $H_{\mathrm{canon}}$ \eqref{Hcanon} by minimizing \eqref{eq:cost_canonHam}, as described in Section~\ref{sec:Honly}. 

We generate $M = 6000$ uniform random vector field samples $(X^{(\mu)}, F^{(\mu)})_{\mu=1}^M$ in the box $(q_1,q_2,p_1,p_2) \in [-1.8,1.8]\times [-2,2]\times [-1.2,1.2]\times [-1.2,1.2]$. We generate $K_{\text{pool}}$ random feature maps and reduce this to $K$ using OMP. 

\begin{figure}[!ht]
\centering
\includegraphics[width=0.8\textwidth]{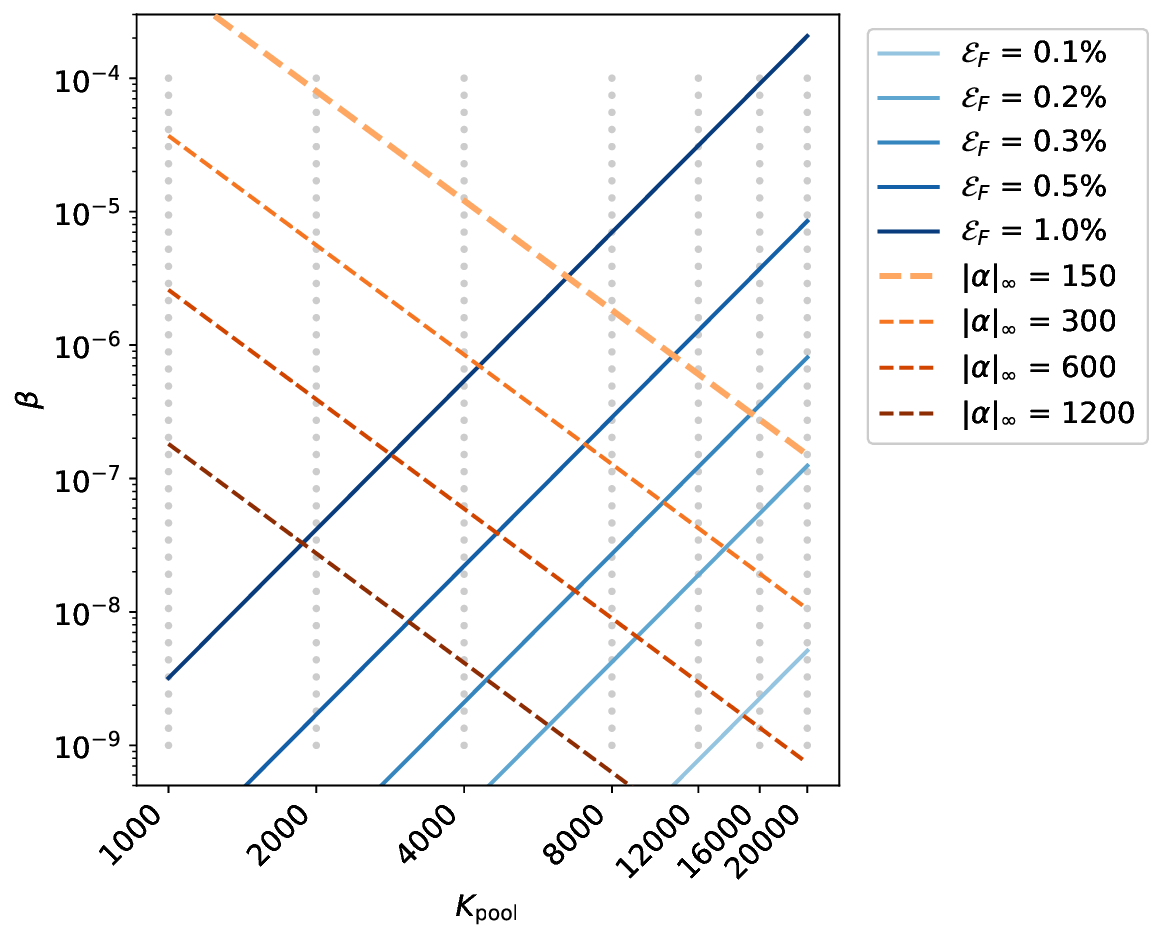}
\caption{Effect of the hyperparameters $K_{\text{pool}}$ and $\beta$ on the accuracy of the learned Hamiltonian of the L86 model~\eqref{L86canon}. We show best fit isolines of the vector field error (blue) and of the maximal learned weights $|\alphaB|_\infty$ for various choices of RFM pool size $K_{\text{pool}}$ and ridge parameter $\beta$. The lines show best fits obtained from results of learned models shown as points in vertical columns.} 
\label{fig:isolines}
\end{figure}

Our method contains two hyperparameters which require tuning: $K_{\text{pool}}$ and $\beta$. We now investigate to what extent they influence the performance. We conducted a parameter study, choosing 
\[
    K_\text{pool}\in \{ 1000,\, 2000,\, 4000, \, 8000, \, 12000, \, 16000, \, 20000 \},\qquad 
    \log_{10} \beta \in [-9,-4],
\]

and compared the relative error between the true and the learned vector fields 
\begin{equation}\label{vferr}
    \mathcal{E}_F = \frac   {\| \J \Win^\top {\rm{diag}}(\alphaB) \PsiB - \F\|_F}{\|\F\|_F}
\end{equation}
and the learned coefficent of largest magnitude $|\alphaB|_\infty$.  The latter quantity affects the stiffness of time integration and is therefore of significance for accurate computation of Poincar\'{e} sections. For each value of $K_\text{pool}$ we learn the model for several values of $\beta$ and record the vector-field error $\mathcal{E}_F$ and $|\alphaB|_\infty$. Figure \ref{fig:isolines} shows best-fit isolines of the vector-field error and of $|\alphaB|_\infty$. The figure indicates how to choose the initial pool size $K_{\text{pool}}$ and ridge parameter $\beta$ to learn the vector field to a desired accuracy, but also provides a rough indication of how the choice affects the stiffness of the resulting learned vector field. 

To reduce the number of feature maps we apply orthogonal matching pursuit (OMP) using the batch-OMP variant of \cite{RuZiEl08}. Use of OMP gives good approximation error with a significantly reduced feature pool as the following table shows (fixed $\beta=3\times 10^{-7}$, least squares fit $K \approx 64.1 K_{\text{pool}}^{0.254}$):

\begin{center}\begin{tabular}{l|rrrrr}
$K_{\text{pool}}$ & 4000 & 8000 & 12000 & 16000 & 20000 \\
\hline
$K$ & 524 & 635 & 686 & 763 & 783
\end{tabular}
\end{center}

\begin{figure}[!ht]
\begin{centering}
\includegraphics[width=0.8\textwidth]{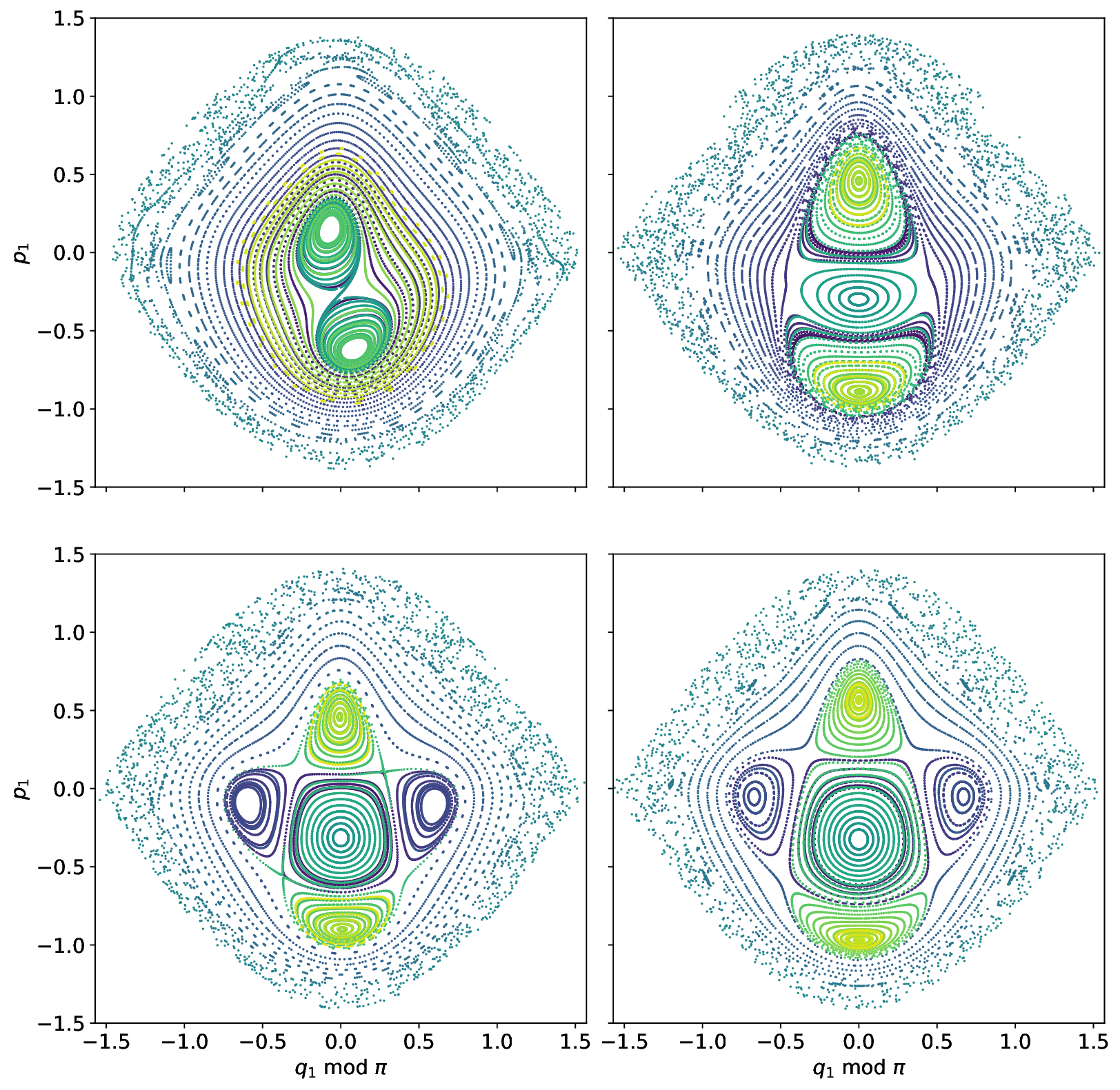}
\caption{Convergence of Poincar\'{e} sections of the L86 model~\eqref{L86canon} in canonical coordinates with learned Hamiltonian, for a range of increasing $K_\text{pool}$, decreasing $\beta$ along the isoline $|\alphaB|_\infty=150$ of Figure \ref{fig:isolines}. $(K_\text{pool},\beta)$ values for upper left: $(2000,8.18\times 10^{-5})$, upper right: $(4000,1.24 \times 10^{-5})$, lower left: $(8000,1.88\times 10^{-6})$, lower right: $(16000,2.86\times 10^{-7})$. We used $M=6000$ training samples and integration stepsize $\Delta t=0.00125$.  
\label{fig:Kpool-sweep}}
\end{centering}
\end{figure}

Using the exponential fit for fixed $|\alphaB|_\infty = 150$ (i.e.~along a line of approximately constant stiffness) we generated four Poincare sections using RFM pool sizes $K_{\text{pool}} = \{ 2000,\, 4000,\, 8000,\, 16000 \}$ (actual selected pool sizes $K=\{ 424, 524, 635, 763\}$, upper orange line in Figure \ref{fig:isolines}) and (fixed) stepsize $\Delta t = 0.00125$. Figure \ref{fig:Kpool-sweep} shows convergence of the Poincare sections. As the RFM pool increases and $\beta$ decreases, the Poincare sections converge to the correct topology of invariant tori shown in Figure~\ref{fig:exactPoincare}.  The vector field errors \eqref{vferr} for these four cases are, respectively, $0.0428$, $0.0181$, $0.0070$, $0.0021$.

\subsection{Inferring the Hamiltonian and the symplectic structure}

In Section \ref{sec:nondegenHam} we proposed a method for jointly inferring the symplectic structure matrix $J$ and the Hamiltonian $H$ for a system for which a canonical relationship among the variables may not be known. We consider again the canonical L86 model~\eqref{L86canon}. 

To test the method, we generate a random matrix $L$ with condition number $2.6667$.
We generate training data using two different approaches. 
First, given $M$ uniform random samples 
$(q_1,q_2,p_1,p_2) \in [-1.8,1.8]\times [-2,2]\times [-1.2,1.2]\times [-1.4,1.4]$, we generate transformed data 
$\tilde{X}^{(\mu)} = L X^{(\mu)}$ and $\tilde{F}^{(\mu)} = L F^{(\mu)}$, $\mu=1,\dots,M$ and use this data to generate an RFM pool using the one shot hit-and-run algorithm~\cite{MaGoNi26}.  

If $K \ge dM$, the matrix of the eigenvalue problem \eqref{Jeig} is identically zero for $\beta=0$. Even for small $\beta>0$ the eigenproblem is overwhelmed by noise without sufficient training data. Because we solve for $\JBhat$ before applying OMP, we require $M > K_\text{pool}/d$. However, a further increase in $M$ has only a mild effect, in our experience. For $K_\text{pool} = 12000$ and $M=3600$, relative error in $\JBhat$ is $10^{-9}$. Error in $\JBhat$ decreases very rapidly, exponentially fit at a rate of about $10^{-4} K_{\text{pool}}$.

\begin{figure}[!ht]
\begin{centering}
\includegraphics[width=0.8\textwidth]{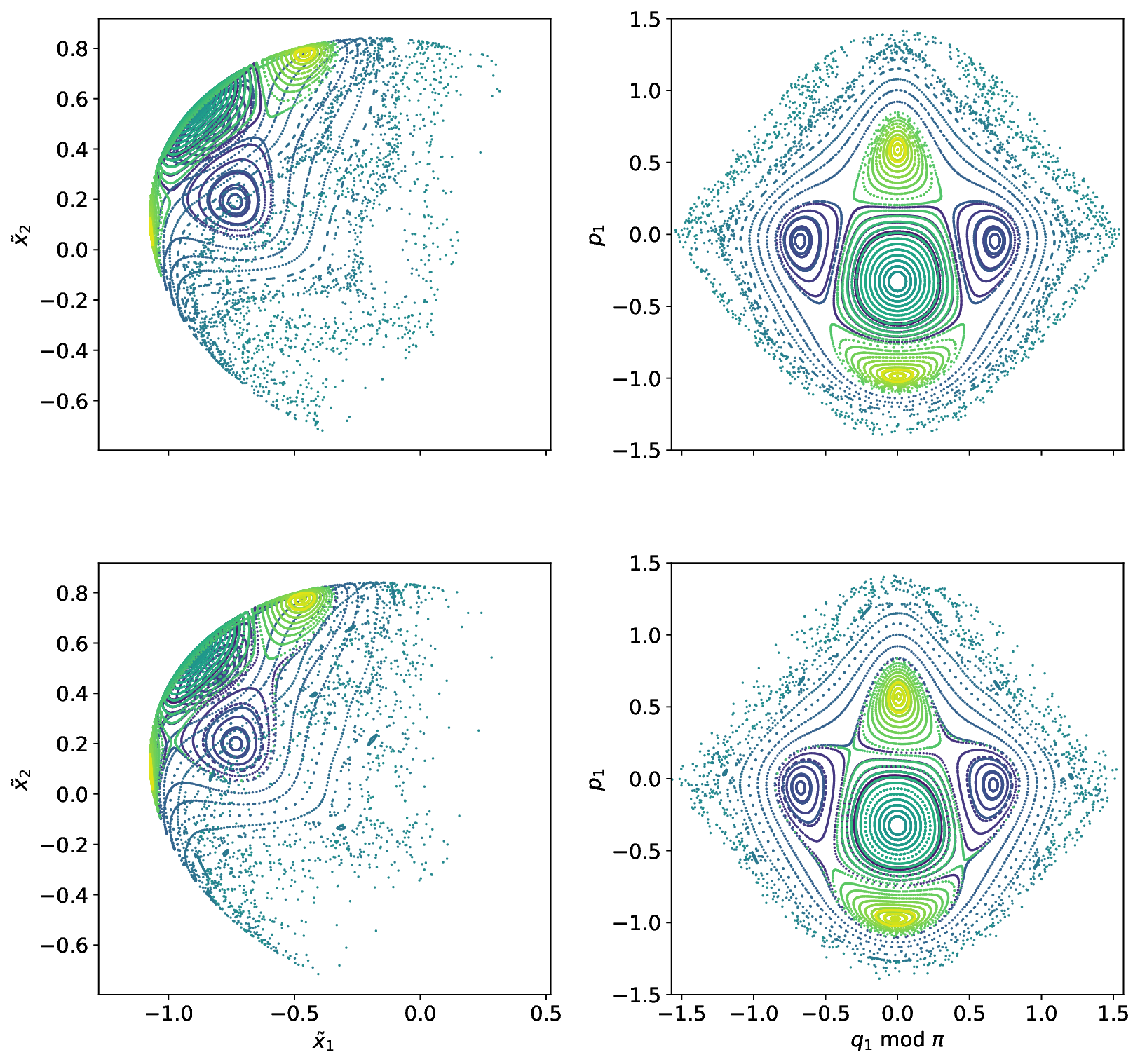}
\caption{Poincar\'{e} sections of the L86 model~\eqref{L86canon} in transformed and in canonical coordinates. The left columns shows the Poincar\'{e} section in the transformed coordinates of the training data and learned system. The right column shows the same Poincar\'{e} sections after applying the inverse transformation. The top row shows simulations of the exact Hamiltonian using a Hamiltonian (Strang) splitting with stepsize $\Delta t = 0.0025$. The lower row shows the inferred Hamiltonian system using the splitting integrator \eqref{Strang} for $K_\text{pool} = 20000$ random feature maps, $\beta = 5 \times 10^{-8}$, $M = 6000$ training samples,  and stepsize $\Delta t = 0.0003125$.
\label{fig:transformed}}
\end{centering}
\end{figure}

Figure \ref{fig:transformed} shows Poincar\'{e} sections of the exact system (top row) and the system with learned (transformed) symplectic structure $\JBhat$ (bottom row), where we used $K_{\text{pool}} = 20000$ and $\beta = 5\times10^{-8}$. On the left, the Poincare section is shown in the transformed coordinates of the training samples and learned dynamics and on the right after inverse transformation to canonical coordinates. The invariant tori are a slightly skewed and the topology is also not completely converged. 

\begin{figure}[!ht]
\begin{centering}
\includegraphics[width=0.8\textwidth]{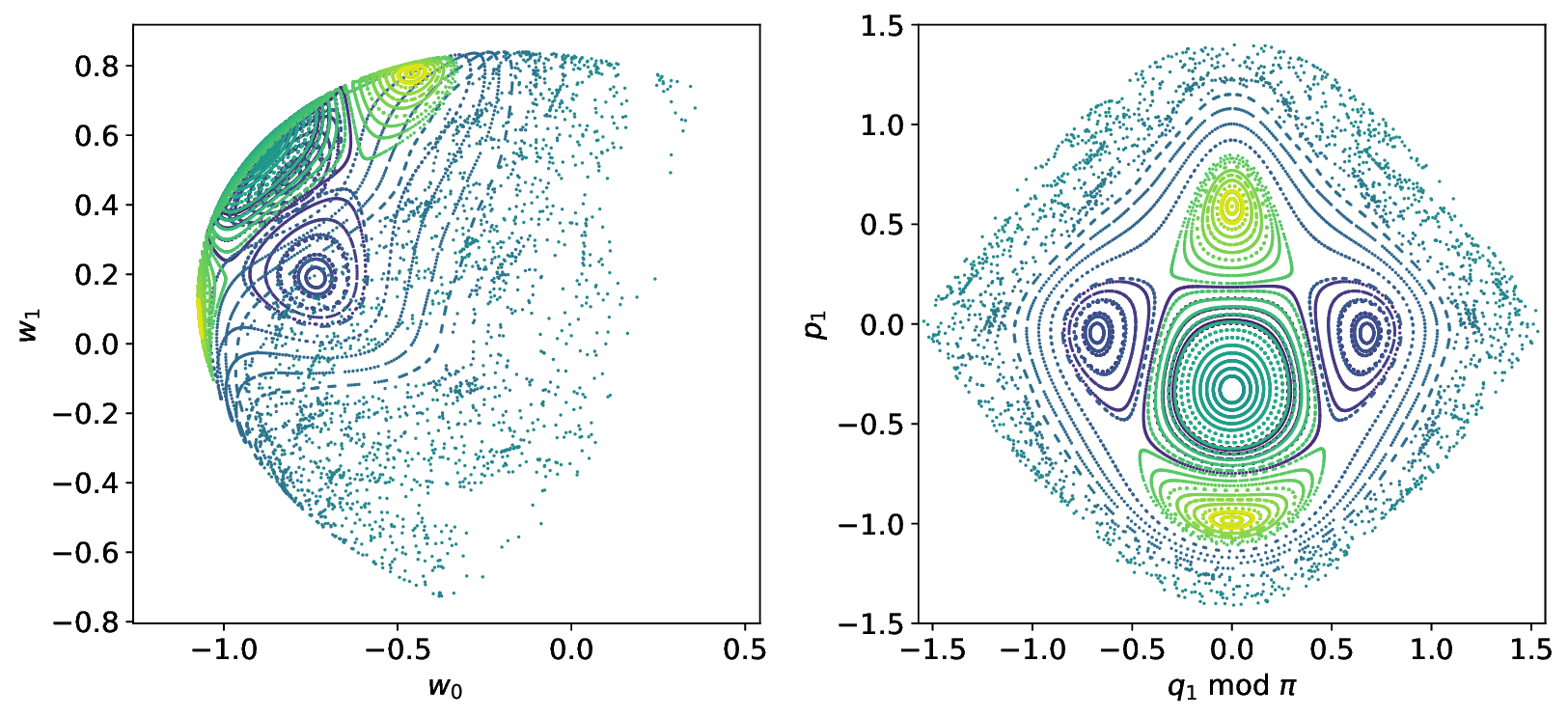}
\caption{Poincar\'{e} sections of the learned L86 model~\eqref{L86canon} in transformed and in canonical coordinates. The training data was generated by running transformed exact model orbits with initial conditions lying on the Poincare section. We used $K_\text{pool} = 12000$ random feature maps, $\beta=5\times 10^{-8}$, $M = 6000$ training samples and stepsize $\Delta t = 0.00125$.
\label{fig:transformed2}}
\end{centering}
\end{figure}

In Figure \ref{fig:transformed2} we repeat the above experiment using alternative training data generated by computing orbits of the transformed exact system using the initial conditions that produce the Poincare section.  In this way, the training data samples precisely the part of phase space that includes the Poincare section. In this case the section is very well reproduced, despite a smaller initial feature pool with $K_\text{pool} = 12000$ and a larger stepsize $\Delta t = 0.00125$.

\subsection{Inferring the Hamiltonian and an affine Poisson structure}

In this section we test the full methodology of Section~\ref{sec:nllsqr} to infer the Hamiltonian and structure constants $C_{ij}^k$ and $J_{ij}$ of the affine Poisson description of the L86 model~\eqref{L86}. We choose an initial RFM pool of $K_{\text{pool}} = 500$ (no OMP), ridge parameter $\beta = 10^{-7}$, $M=4000$ training points drawn uniformly from a rectangular box $x\in [-1.5,1.5]^5$. The nonlinear least squares of Section \ref{sec:nllsqr} converges to relative vector field error $0.00036$. The structure coefficients are determined to relative errors of $2.2\times 10^{-4}$ for $C$ and to $8\times 10^{-5}$ for $J$.

\begin{figure}[!ht]
\begin{centering}
\includegraphics[width=0.8\textwidth]{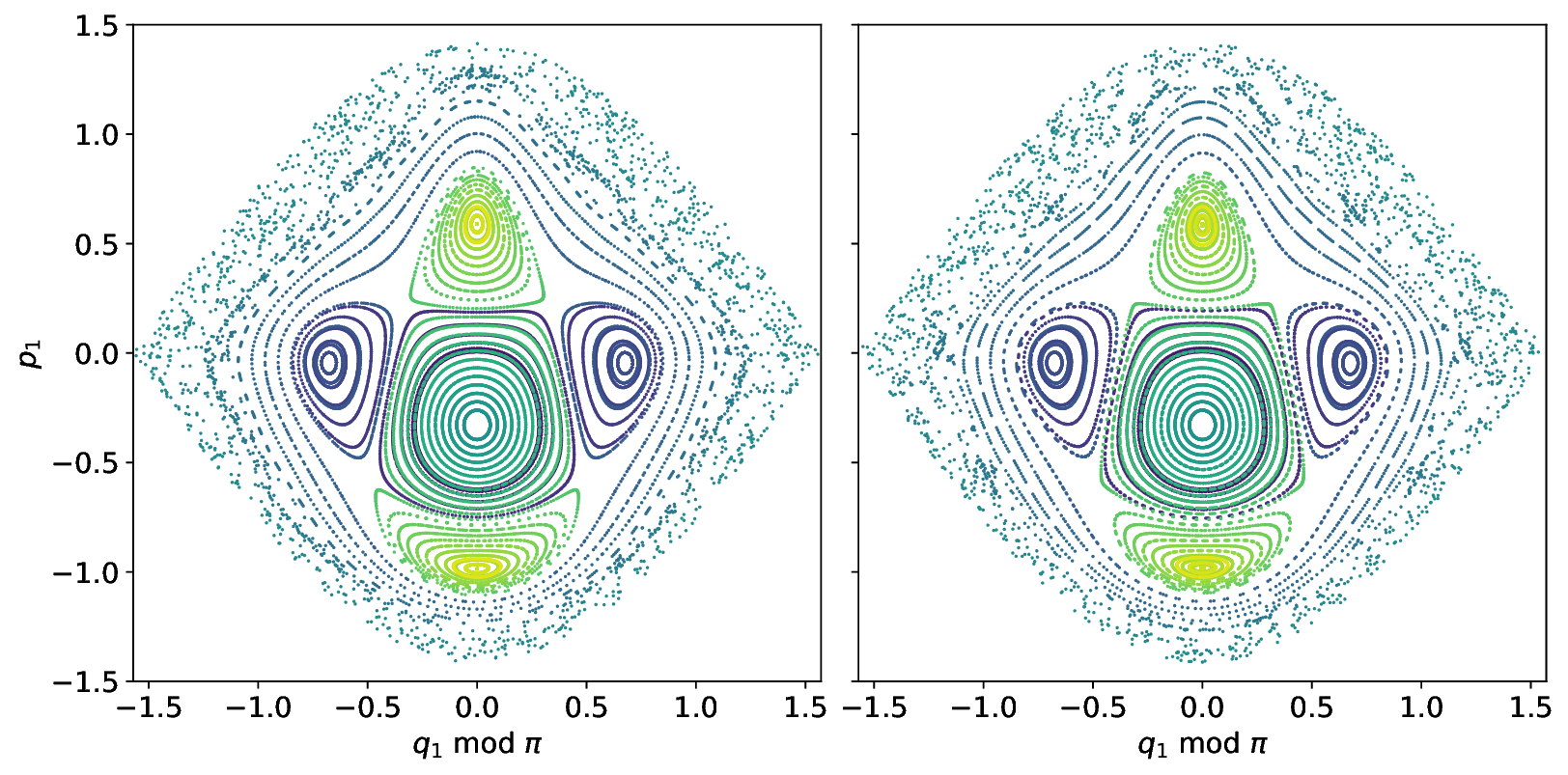}
\caption{Poincar\'{e} section of the learned L86 model~\eqref{L86} on $\R^5$ with affine Poisson structure \eqref{L86affineB} for the exact system (left) and the learned model (right). For the true model we used a stepsize $\Delta t = 0.0025$. For the learned model we used  $K = 500$ features (no OMP), $\beta = 10^{-7}$, $M = 4000$ training samples and a stepsize $\Delta t=0.02$. Both systems were integrated using a symmetric Poisson splitting integrator.
\label{fig:affine}}
\end{centering}
\end{figure}

We integrate the learned model using the splitting method of Section \ref{sec:splitting}. We constructed the Poincar\'{e} sections of the true system \eqref{L86} using stepsize $\Delta t = 0.0025$ and the learned model using a stepsize $\Delta t = 0.02$ in the $\R^5$ description.  Figure \ref{fig:affine} shows the comparison, which is already very close at this relatively low resolution fit.

\begin{figure}[!ht]
\begin{centering}
\includegraphics[width=0.8\textwidth]{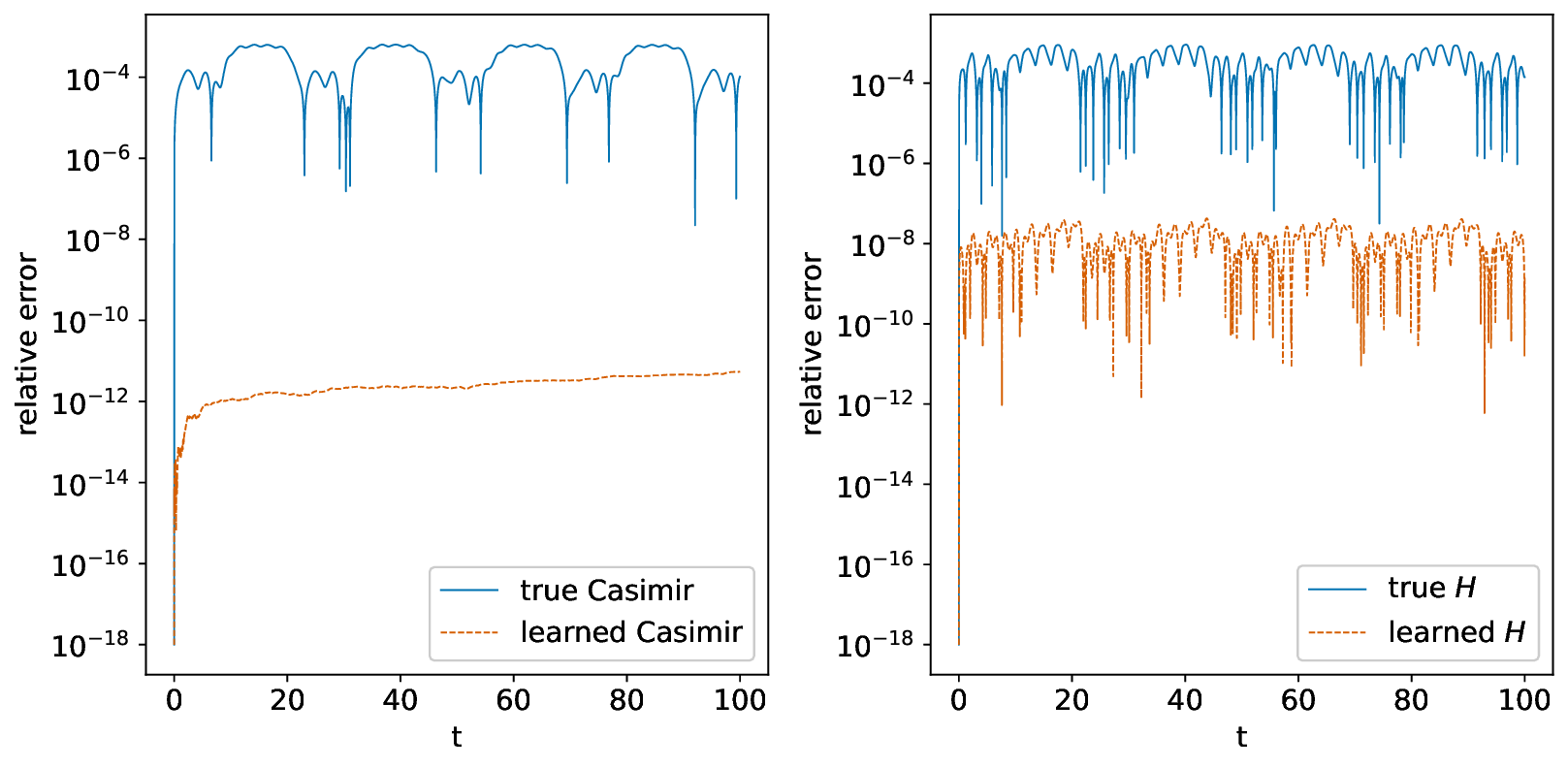}
\caption{Relative errors in conservation laws for the learned L86 model with affine Poisson structure: (left) the error in the true system's Casimir \eqref{L86cas} and the exact Casimir of the learned system; (right) the error in the true system's Hamiltonian \eqref{L86ham} and the learned Hamiltonian.\label{fig:L86conserve}}
\end{centering}
\end{figure}

The Casimir of the learned model was found to be slightly perturbed with respect to the true $Z$. It is  quadratic in $x$, but with a small dependence on all components and it contains an additional small linear term.  Figure \ref{fig:L86conserve} shows the error in the Casimir and Hamiltonian with respect to the first integrals of the true system and those of the learned system.  The relative errors in the conserved quantities of the true system are both about $10^{-3}$.  The relative errors of the exact conserved quantities of the learned system are significantly smaller:  about $10^{-7}$ for the Hamiltonian and about $10^{-12}$ for the Casimir. The effects of many rounding errors is evident in the slow error growth of the Casimir $Z$.

\subsection{Inferring an affine Poisson structure for a quadratic Hamiltonian}

Many physical Hamiltonian systems have a quadratic Hamiltonian. As discussed in Section \ref{sec:Ham}, affine Poisson (Lie-Poisson) systems remain affine Poisson (Lie-Poisson) under a linear change of coordinates. Consequently, if the Hamiltonian of the true dynamics is quadratic with $H(x) = \frac{1}{2} x^\top M x + c^\top x$ (with $c=0$ in the Lie-Poisson case), where $M$ is symmetric positive definite, then there exists a linear change of variables $\tilde {x} = L x + y$ (with $y=0$ in the Lie-Poisson case) such that the Hamiltonian becomes isotropic in the sense that $H(\tilde{x}) = \frac{1}{2} \| \tilde{x} \|^2$. This transformation implies that without loss of generality any system with a quadratic Hamiltonian can be assumed to be of this form. Hence, the system identification problem reduces to one where only the symplectic structure matrix or the Poisson tensor needs be learned from data. This assumption significantly reduces the computational complexity and allows for higher accuracy. We illustrate this for the affine representation of the L86 model \eqref{L86}.  

\begin{figure}[!ht]
\begin{centering}
\includegraphics[width=0.5\textwidth]{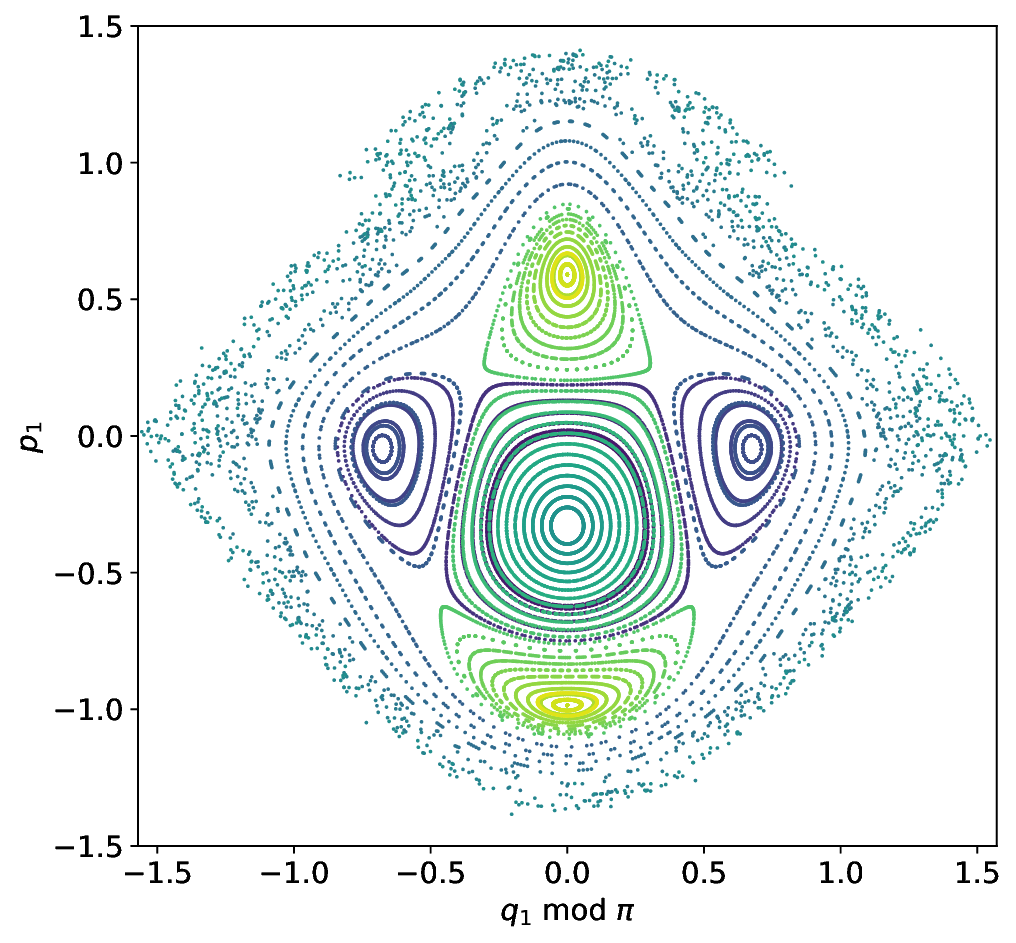}
\caption{Poincar\'{e} section of the learned L86 model~\eqref{L86} for which the affine Poisson structure is inferred assuming a quadratic Hamiltonian. We used $M = 2000$ data points and a stepsize $\Delta t = 0.01$.
\label{fig5:quadraticH}}
\end{centering}
\end{figure}

Assuming Hamiltonian $H = \frac{1}{2} \|x\|^2$, we can replace $\nabla \Hhat(\X^{(\mu)})$ with $\X^{(\mu)}$ in the loss function \eqref{L3}, thereby eliminating the feature maps and $\alphaB$ altogether. We solve directly for the structure parameters $\thB = (\Chat_{ij}^k,\Jhat_{ij})$ using the continuation method of Section \ref{sec:nllsqr}. We implemented this for the L86 model~\eqref{L86}, drawing $M=2000$ vector field data uniformly on $x \in [-1.5,1.5]^5$.  The transformed (due to isotropization) structure constants are exact to machine precision relative error.  Figure \ref{fig5:quadraticH} shows the Poincare section. 

\subsection{Inferring the Hamiltonian and the Lie-Poisson structure}

We now test the methodology of Section \ref{sec:LPstruct} on Lie-Poisson systems. As stated in Section \ref{sec:Ham}, the number of Jacobi identity constraints exceeds the number of structure coefficients for dimension $d>5$. In this section we present numerical results for the Kirchhoff-type model~\eqref{kirchhoff} for which $d=6$. For the numerical experiment we chose $a = (0.5452, 0.8285, 0.5628)$, $b=(-0.0756, -0.2957, 0.2476)$, and $d = (0.6759, 1.1507, 0.7016)$, and initial condition
\[
    \Omega(0) = (0.5, 0.3, -0.4), \qquad \Gamma(0) = (1.0, 0.2, -0.3).
\]

We generated $M=6000$ phase space samples uniformly in the hypercube $x = (\Omega, \Gamma) \in [ -1.5, 1.5 ]^6$ and learned the Lie algebra structure constants $\Chat_{ij}^k$ and Hamiltonian $\Hhat$ using $K=500$ random feature maps (without OMP) with ridge parameter $\beta = 10^{-7}$.
The learned structure constants matched the true structure constants with a maximum entrywise error of 0.000132. 

\begin{figure}[!ht]
\begin{centering}
\includegraphics[width=0.8\textwidth]{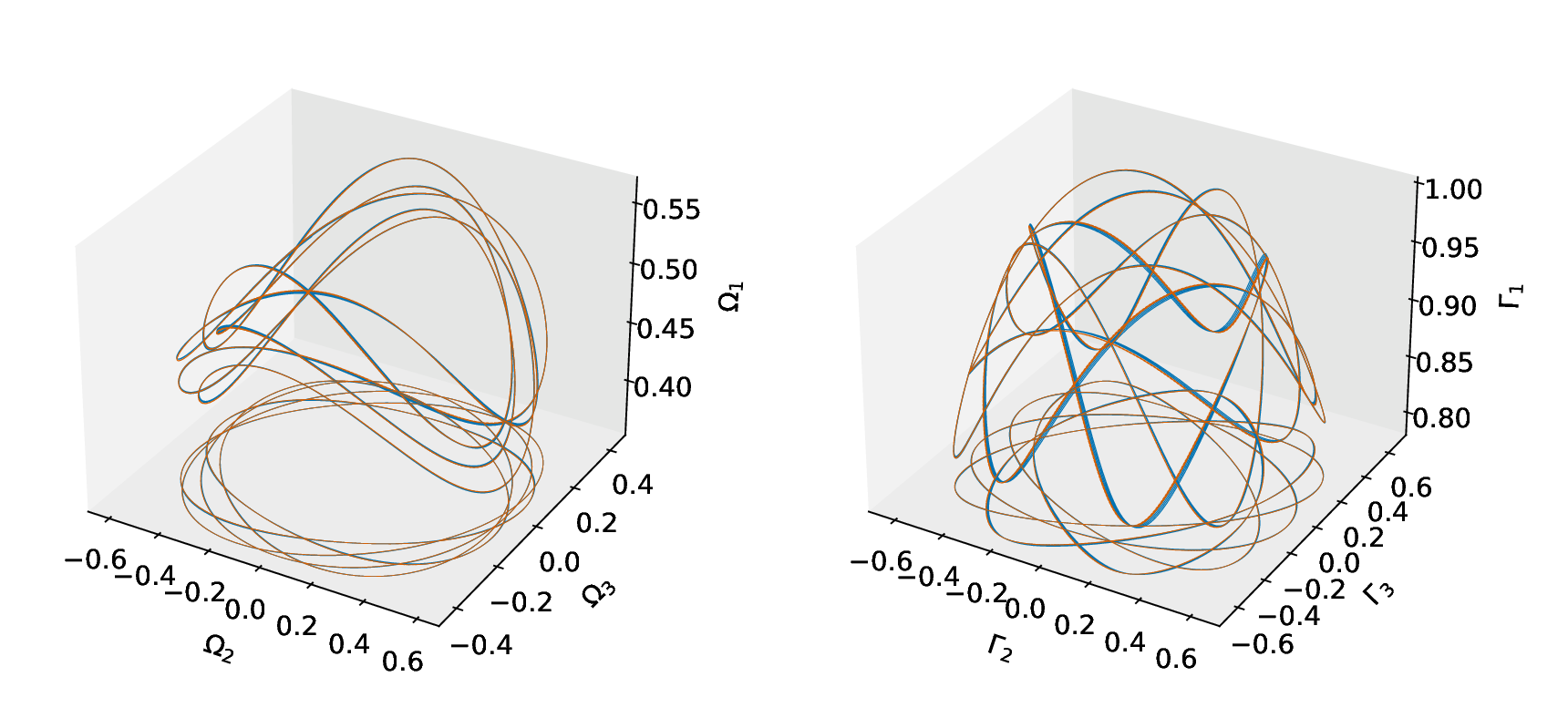}
\caption{Trajectories of the true (blue) and learned (orange) dynamics of $\Omega(t)$ (left) and $\Gamma(t)$ (right). We used $K=500$ features (no OMP), $\beta=10^{-7}$, $M=6000$ training samples, integration stepsize $\Delta t = 0.01$ and orbit length $T=300$.
\label{fig:se3}}
\end{centering}
\end{figure}

Figure \ref{fig:se3} illustrates the computed orbits of the true and of the learned systems, integrated with the Poisson splitting method \eqref{Strang} described in Section \ref{sec:splitting} with stepsize $\Delta t=0.01$ up to time $T=300$.  The trajectories show good agreement.

The learned system admits two Casimirs, but these may differ significantly from those of the true system due to the ambiguity of multiple Casimirs. 
That is, since the Poisson structure of the learned system is an approximation of the true Poisson structure, its Casimirs could be approximations of combined functions of the true system's Casimirs:
\[
    \tilde{Z}_1 \approx g_1 (Z_1,Z_2), \quad
    \tilde{Z}_2 \approx g_2 (Z_1,Z_2).
\]

\begin{figure}[!ht]
\begin{centering}
\includegraphics[width=0.8\textwidth]{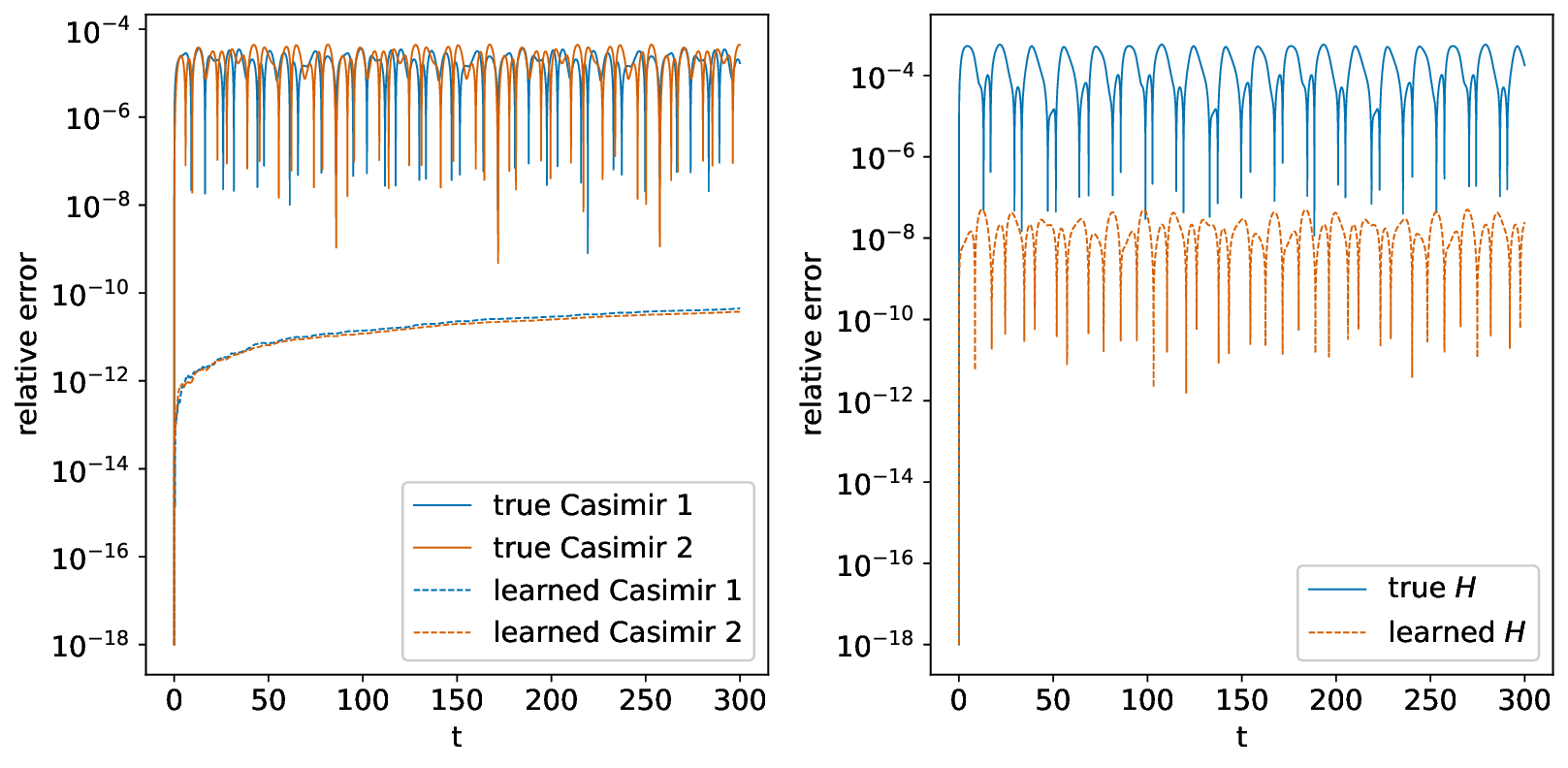}
\caption{Error in conservation of Casimirs (left) and Hamiltonian (right) for the Kirchhoff model~\eqref{kirchhoff} by the learned model. Due to the ambiguity between the Hamiltonian and Casimirs, the learned model's exact first integrals differ from those of the true model. In the figures we indicate conservation errors with respect to the true model's first integrals and those of the learned model.  Integration stepsize stepsize is $\Delta t=0.01$.\label{fig:se3conserve}}
\end{centering}
\end{figure}

Figure \ref{fig:se3conserve} shows errors in the Casimirs and in the Hamiltonian of the true system, evaluated along the learned system's trajectory, as well as the error in conservation of the learned Hamiltonian and Casimirs of the learned Poisson tensor.  The learned system's Casimirs are conserved to relative error less than $10^{-10}$ over the whole time window. The observed drift grows at order $t^{1/2}$ consistent with random rounding errors due to the large number of terms in the splitting; we employ $2K = 1000$ matrix exponential computations per time step.

\section{Discussion and outlook}
\label{sec:disc} 

We have developed a data-driven methodology for system identification of general symplectic Hamiltonian systems and for Hamiltonian systems with an affine Poisson structure, including Lie-Poisson systems. We achieved this by expressing the Hamiltonian as a random feature map. This allows for a computationally efficient estimation with an explicit solution of the least-squares regression problem in the case when the Poisson tensor is given. To learn both the Hamiltonian and the Poisson tensor we employed a bilevel optimization. The nonsmooth nature of the algebraic variety defined by the Jacobi identity for affine Poisson systems posed a particular challenge that we solved by an gradual enforcement of the constraint. We provided numerical results that show that the methodology successfully recovers the true dynamics even in systems with dimension larger than five for which the number of constraints exceeds the number of unknowns. We remark that we have done additional successful experiments with non-physical randomly generated Lie-Poisson systems with dimensions up to $\R^9$. 

The methodology as introduced here relies on several assumptions, that we now discuss. We assume here that we have access to the whole state variables $x\in \R^{2d}$ so in particular, we assume that we observe the conjugate momenta in the canonical setting. This may not be the case in a real-world situation. 

Further, we have considered here the case where we are given samples of the state space variables and the vector field, i.e.~the time derivative of the vector field. In practice one typically faces the situation where one has access to noisy measurements of the state space variables only. In this case, one would need to estimate the vector field from the state space variables. The presence of noise prohibits any estimation of the vector field via finite-time discretizations, even for small sampling times. One may either employ denoising via splines \cite{KaKuBr20}, control the noise via data assimilation \cite{GottwaldReich21a} or by Gaussian process regression~\cite{HuOrYi25a}. The important case of partial observations is also not considered here and may be tackled via delay embedding \cite{GottwaldReich21}. 

We chose here RFMs to express the Hamiltonian. One may instead express the Hamiltonian as a linear combination of functions from a prescribed library. This opens up the way to perform SINDy~\cite{BrPrKu16} and achieve a sparse identification of Lie-Poisson systems. 

An interesting application of our framework would be to actually discover if a dynamical system is Lie-Poisson. For the heavy top, for example, we could take measurements in Cartesian coordinates. The associated Hamiltonian system in these coordinates is then not Lie-Poisson. However, by allowing for a diffeomorphic transformation $z=T(x)$ we can test if the dynamical system allows for a Lie-Poisson structure. We can represent the transformation $T$ by a random feature map and impose the constraints for the structure matrix to be Lie-Poisson as described in Section~\ref{sec:LPstruct}. 
However, the transformation is now entering the nonlinear activation function and the loss function will therefore not be a quadratic function of the coefficients of the random feature map for the transformation. 

{\bf Repository and AI statement.} 
Codes to generate all data and figures in this article are available at \url{https://doi.org/10.5281/zenodo.24247926}.
Sonnet 5.5 by Anthropic was used to port the authors' original Matlab codes to Python and C, as well as to interface the execution of numerical experiments. The authors take all responsibility for correctness of the reported results.


\bibliographystyle{abbrv}
\bibliography{FrGo26}
 
\end{document}